\documentclass[11pt]{amsart}
\usepackage{latexsym,amssymb,amsmath}
\usepackage{youngtab}
\input epsf
\def\boxit#1{\vbox{\hrule height1pt\hbox{\vrule width1pt\kern3pt
  \vbox{\kern3pt#1\kern3pt}\kern3pt\vrule width1pt}\hrule height1pt}}

\def\trank{\text{rank}}

\def\BC{\mathbb C}\def\BF{\mathbb F}\def\BO{\mathbb O}\def\BS{\mathbb S}
\def\BA{\mathbb A}\def\BQ{\mathbb Q}
\def\BP{\mathbb P}
\def\pp#1{\mathbb P^{#1}}
\def\fa{\mathfrak a}\def\fr{\mathfrak r}\def\fz{\mathfrak z}
\def\fb{\mathfrak b}\def\fo{\mathfrak o}\def\fgl{\mathfrak g\mathfrak l}

\def\fd{\mathfrak d}

\def\ppp{{\mathbb P}}
\def\pp#1{{\mathbb P}^{#1}}
\def\tdim{\rm dim}
\def\hd{,...,}
\def\ww{\wedge}
\def\upperp{{}^\perp}

\def\na{n+a}

\def\inv{{}^{-1}}
\def\com{{\rm com}}
\def\trace{{\rm trace}}
\def\cI{{\mathcal I}}\def\cB{{\mathcal B}}\def\cA{{\mathcal A}}
\def\cJ{{\mathcal J}}
\def\cE{{\mathcal E}}\def\cT{{\mathcal T}}
\def\cF{{\mathcal F}}
\def\cG{{\mathcal G}}
\def\cR{{\mathcal R}}
\def\cS{{\mathcal S}}
\def\cL{{\mathcal L}}
\def\cW{{\mathcal W}}
\def\cO{{\mathcal O}}
\def\CC{\mathbb C}
\def\RR{\mathbb R}
\def\HH{\mathbb H}
\def\AA{{\mathbb A}}
\def\BB{{\mathbb B}}
\def\OO{\mathbb O}
\def\LG{\mathbb {LG}}
\def\LF{{\mathbb {LF}}}
\def\ZZ{\mathbb Z}
\def\SS{\mathbb S}
\def\GG{\mathbb G}
\def\11{\mathbf 1}
\def\PP{\mathbb P}
\def\QQ{\mathbb Q}
\def\FF{\mathbb F}
\def\JA{{\mathcal J}_3(\AA)}
\def\JB{{\mathcal J}_3(\BB)}
\def\ZA{{\mathcal Z}_2(\AA)}
\def\fh{{\mathfrak h}}
\def\fs{{\mathfrak s}}
\def\fsl{{\mathfrak {sl}}}
\def\fsp{{\mathfrak {sp}}}
\def\fspin{{\mathfrak {spin}}}
\def\fso{{\mathfrak {so}}}
\def\fe{{\mathfrak e}}

\def\ff{{\mathfrak f}}

\def\fz{{\mathfrak z}}
\def\ffi{{\mathfrak i}}
\def\fg{{\mathfrak g}}
\def\fn{{\mathfrak n}}
\def\fp{{\mathfrak p}}
\def\fk{{\mathfrak k}}
\def\ft{{\mathfrak t}}
\def\fl{{\mathfrak l}}
\def\l{\lambda}
\def\a{\alpha}
\def\ta{\tilde{\alpha}}
\def\o{\omega}
\def\oo{\Omega}
\def\O{\Omega}
\def\b{\beta}
\def\g{\gamma}
\def\s{\sigma}
\def\k{\kappa}

\def\d{\delta}
\def\th{\theta}
\def\m{\mu}
\def\up#1{{}^{({#1})}}
\def\e{\varepsilon}
\def\ot{{\mathord{\,\otimes }\,}}
\def\op{{\mathord{\,\oplus }\,}}
\def\otc{{\mathord{\otimes\cdots\otimes}\;}}
\def\pc{{\mathord{+\cdots +}}}
\def\lra{{\mathord{\;\longrightarrow\;}}}
\def\ra{{\mathord{\;\rightarrow\;}}}
\def\da{{\mathord{\downarrow}}}
\def\we{{\mathord{{\scriptstyle \wedge}}}}
\def\JA{{\mathcal J}_3(\AA)}
\def\JB{{\mathcal J}_3(\BB)}
\def\tr{{\rm trace}\;}
\def\dim{{\rm dim}\;}
\def\La#1{\Lambda^{#1}}

\newcommand{\norm}[1]{\lVert#1\rVert}

\newtheorem{theo}{Theorem}[section]
\newtheorem{coro}[theo]{Corollary}

\newtheorem{prop}[theo]{Proposition}

\newtheorem{theorem}{Theorem}[section]
\newtheorem{proposition}[theorem]{Proposition}

\theoremstyle{definition}

\newtheorem{question}[theorem]{Question}
\newtheorem{example}[theorem]{Example}

\theoremstyle{remark}
\newtheorem{remark}[theorem]{Remark}

\begin{document}

\title{Representation theory and projective geometry}
\author{J.M. Landsberg and L. Manivel}
\date{July 2002}
\maketitle

{\tableofcontents}

\section{Overview}

This article has two purposes.  The first is to provide an elementary
introduction to papers \cite{LM0,LM1,LMmagic,LMtrial,LMseries},
 related works and their historical context.
  Sections 2.1--2.3, 2.5--2.6, 3.1--3.5, 4.1--4.3 should be accessible to a
general audience of mathematicians.   The second is to provide generalizations,
new perspectives, and complements to results in these papers, i.e., things we
thought of after the papers were published. In particular, we mention  2.5, 2.7, 3.4 and 
3.6--3.9. Each section begins with a description of its
contents.
 
\medskip

Simply put, our goals are to use geometry to solve questions in representation
theory and to use representation theory to solve questions in geometry.

On the geometric side, the objects of interest are homogeneous
varieties $X=G/P\subset \ppp V$. Here  $G$ is a complex semi-simple
Lie group, $V$ is an irreducible $G$-module,   $X$ is the unique closed
orbit (projectivization of the orbit of a highest weight vector)
and $P$ is the stabilizer of a point.
For example, let $W=\BC^m$ and let $X=G(k,W)\subset\ppp (\La k W)$, the
Grassmannian of $k$-planes through the origin in $W$.
Here $G/P=SL(W)/P$, $V=\La kW$.
We are more generally interested in the geometry of orbit closures in
$\ppp V$. 

\smallskip

Basic questions one can ask about a variety are its dimension, its
degree, and more generally its Hilbert function. We may ask
the same questions for varieties associated to $X$. For example the
degrees of the dual varieties of Grassmannians are still unknown
in general (see \cite{lascoux, dw,tel}).  

Other types of problems include recognition questions. For example,
given a variety with certain geometric properties, are those properties
enough to characterize the variety? For an example of this, see
Zak's theorem on Severi varieties below in \S 2.6. For another example,
 Hwang and Mok characterize rational homogeneous
varieties by the variety of tangent  directions to minimal degree rational
curves passing through a general point, see
\cite{hwangmokpapers} for an overview.  Some results along those
lines are described in \S 2.5 below.
We also mention the LeBrun-Salamon conjecture \cite{LeBrun, LS}
which states that any Fano variety equipped with a holomorphic contact
structure must be the closed orbit in the adjoint representation
$X_{ad}\subset\ppp\fg$ for a complex simple Lie algebra. In this context
also see \cite{bea, kebI}. 

\medskip

On the representation theory side, the basic objects are $\fg$, a 
complex semi-simple  Lie algebra and $V$, an irreducible
$\fg$-module (e.g., $\fg=\fsl (W)$, $V=\La kW$).
Problems include  the classification of orbit closures
in $\BP V$,  to 
construct explicit models for the  group action,
to geometrically interpret the decomposition of
$V^{\ot k}$ into irreducible
$\fg$-modules.
  We discuss
these classical questions below, primarily for algebras
occuring in \lq\lq series\rq\rq .

 Vassiliev theory points to the need for defining objects
more general  than Lie algebras. We have nothing to
add about this subject at the moment, but the results of
\cite{LMtrial, LMseries} were partly inspired by   work of Deligne \cite{del}
and Vogel \cite{vog1,vog2} in this direction.

 For the mystically inclined, there are many strange
 formulas related to the exceptional groups.
We present some such formulas in \S 4.3, 4.6 below.  
  Proctor and Gelfand-Zelevinski  filled in  
\lq\lq holes\rq\rq\ in the classical formulas for the
$\fo\fs\fp _n$ series using the non-reductive
{\it odd symplectic groups}. Our formulas  led  us to exceptional 
analogues of the odd symplectic
groups. These analogues are currently under
investigation (see \cite{LMW}).

When not otherwise specified, we use the ordering of roots
as in \cite{bou}.

We now turn to details. We begin with some observations
that lead to   interesting rational maps of projective spaces.
\bigskip

\section{Construction of complex simple Lie algebras via geometry}

We begin in \S 2.1  with three ingredients that go
into our study: local differential geometry (asymptotic
directions), elementary algebraic
geometry (rational maps of projective space)
and   homogeneous
varieties (the correspondence between rational homogeneous
varieties and marked Dynkin diagrams). We then, in \S 2.2-2.4
describe two algorithms that construct new varieties
from old that lead to new proofs of the classification
of compact Hermitian symmetric spaces and the Cartan-Killing
classification of complex simple Lie algebras.
The proofs are constructive, via explicit rational maps and
in \S 2.5 we describe applications and generalizations
of these maps. Our maps generalize maps used by Zak
in his classification of Severi varieties and in \S 2.6 we
describe his influence on our work. In \S 2.7 we return to
a topic raised in \S 2.1, where we determined the parameter space
of lines through a point of a homogeneous variety $X=G/P$. We explain
{\it Tits' correspondences} which allow one to determine the
parameter space of all lines on $X$ and in fact 
parameter spaces for all $G$-homogeneous
varieties on $X$. We explain how to use Tits correspondences to 
explicitly construct
certain homogeneous vector bundles, an in turn to use the explicit
construction to systematize Kempf's method for desingularizing 
orbit closures.

\subsection{Differential geometry,   algebraic
geometry and representation theory}

\subsubsection{Local differential geometry}

Let $X^n\subset\pp\na$ be an algebraic variety.

\begin{question}How to study the geometry of $X$?\end{question}

One answer: Take out a microscope. Fix a smooth point $x\in X $ and   take
derivatives at $x$. The first derivatives don't yield much. One gets
$\tilde T_xX$, the embedded tangent projective space,
  the union of 
lines (linearly embedded
$\pp 1$'s) having contact to order one with $X$ at $x$. (Contact to
order zero means the line passes through $x$.)

Sometimes we   work with vector spaces so for 
$Y\subset \ppp V$ we let $\hat Y\subset V$ denote the cone over $Y$
and we let $\hat T_xX=\hat{\tilde T}_xX$.
Let $T_xX$ denote the Zariski (intrinsic) tangent space to $X$ at $x$. We have
 $T_xX=\hat x^*\ot \hat T_xX/\hat x$, see \cite{LM0}.

Taking second derivatives one obtains a variety $Y=Y_{X,x}\subset
\ppp T_xX$ whose points are the {\it asymptotic directions}, the
tangent directions to lines having contact to order two with $X$ at $x$.

$$\epsfbox{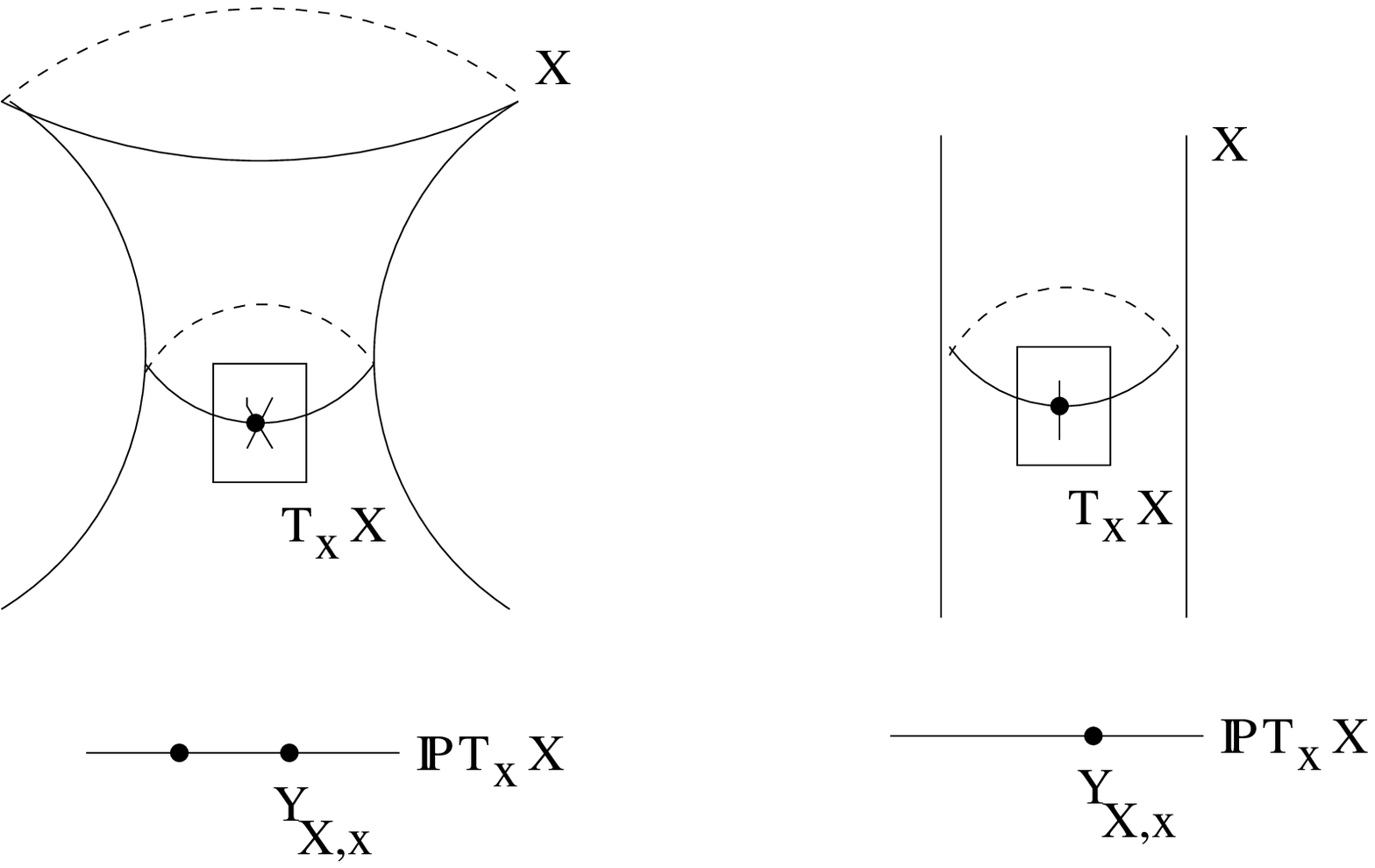} $$

In an attempt to get global
information from the asymptotic directions,
restrict to the case where   $x\in X$ is a {\it general} point

\begin{question} How much of the geometry of $X$ can be recovered from
$Y_{X,x}\subset\ppp T_xX$?\end{question}

Answer: Usually  not much. For example, if $X$ is a smooth hypersurface,
then $Y$ is always a smooth quadric hypersurface, i.e., all smooth hypersurfaces
look the same to second order.

The set of asymptotic directions is the zero set
of a system of quadratic equations generically of dimension equal to
the codimension of $X$ (unless the codimension of $X$ is
large, in which case it is generically the complete system of
quadrics). If $Y$ is sufficiently pathological one might hope to recover
important information about $X$.

Consider the
{\it Segre variety},
 $Seg(\pp k\times \pp l)\subset\ppp (\BC^{k+1}\ot \BC^{l+1})$ of rank one
 matrices in the space of all $(k+1)\times (l+1)$ matrices. A short calculation
 shows   $Y_{Seg,x}=\pp{k-1}\sqcup\pp{l-1}$,
 the disjoint union of a $\pp{k-1}$ with a
 $\pp{l-1}$. (Note that the codimension is sufficiently large here
 that one would expect $Y$ to be empty based on dimension considerations.) 

   Griffiths
and Harris \cite{gh} conjectured that if
$Z \subset\pp 8$ is a variety of dimension $4$
and $z\in Z$ a general point, if $Y_{Z,z}\subset\ppp T_zZ$ is the disjoint union
of two lines, then $Z=Seg(\pp 2\times \pp 2)\subset\pp 8$.

\begin{theorem} \cite{Lrigid}
Let $k,l>1$, let $Z^{k+l}\subset\pp N$ be a variety, 
and let $z\in Z $ be a general point.
If $Y_{Z,z}=\pp{k-1}\sqcup\pp{l-1}$, then $Z=Seg(\pp k\times\pp l)$.
\end{theorem}

Moreover, the analogous rigidity is true for
varieties having the same asymptotic directions
as $X=G(2,m)\subset\ppp(\La 2\BC^m)$,
the Grassmannian of two-planes in $\BC^m$
in its Plucker embedding, with $m\geq 6$,
and the complex Cayley plane $E_6/P_6={\BO\BP^2_\BC}\subset\pp{26}$ 
\cite{Lrigid}. We  recently prove  similar
rigidity  results for homogeneous Legendrian varieties \cite{LMleg}.

The rigidity results hold in the $C^{\infty}$ category if one replaces
the word \lq\lq general\rq\rq\ by the word \lq\lq every\rq\rq .

If one adds the global assumption that $Z$ is smooth, the
results also hold for $X=G(2,5)$ and $X=\BS_5 \subset\pp{15}$,
the spinor variety using different methods.

We are unaware of any second order rigidity results for nonhomogeneous
varieties. Thus, purely from the most naive differential geometry, one
already encounters homogeneous varieties   
as examples of the most rigid projective
varieties.  In fact, so far we have just encountered
the most homogeneous ones, the ones admitting Hermitian symmetric
metrics.
For more on this  and the rigidity of other homogenous varieties
see \S 2.5.

 We define a homogeneous variety $X=G/P\subset \BP V$ to be {\it minuscule}
if $G$ is simple,
$X$  admits a Hermitian symmetric metric and $X$ is in its minimal
homogeneous embedding. $X$ is said to be {\it generalized minuscule}
if it is homogeneously embedded and admits a Hermitian symmetric metric.

\begin{remark} If $X\subset \BP V$ is a variety cut out by quadratic polynomials,
then the asymptotic directions $Y\subset \BP T_xX$ are actually the
tangent directions to lines (linearly embedded $\BP^1$'s) on $X$.
Rational homogeneous varieties are cut out by quadratic equations, in
fact if $X\subset \BP V_{\lambda}$, then the ideal is generated
by $V_{2\lambda}\upperp\subset S^2V_{\lambda}^*$.\end{remark}

\subsubsection{Lie groups and homogeneous varieties}

Let $G$ be a complex (semi)-simple Lie group and let $V$ be an irreducible
$G$-module. Then there exists a unique closed orbit $X=G/P\subset \ppp V$.

\smallskip\noindent {\sl Examples}:

1. $G=SL(n,\BC)$, the group preserving
$det\in \La n\BC^{n*}$, $V=\La k\BC^n$, $X=G(k,n)$, the Grassmannian of
$k$-planes through the origin in $\BC^n$.

2. $G=SO(n,Q)$, the group preserving a nondegenerate $Q\in S^2\BC^{n*}$,
$V=\La k\BC^n$, $X=G_Q(k,n)$, the Grassmannian of $Q$-isotropic
$k$-planes through the origin in $\BC^n$.

3. $G=Sp(n,\o)$, the group preserving a nondegenerate $\o\in \Lambda^2\BC^{n*}$,
$V=\La k\BC^n/(\La{k-2}\BC^n \ww\o)$, $X=G_{\o}(k,n)$, the Grassmannian of $\o$-isotropic
$k$-planes through the origin in $\BC^n$. Here $n$ is usually required
to be even (but see \cite{proc, gz}, and \S 4.3).

\medskip

Since linear algebra is easier than   global
geometry, we work with
$\fg=T_{Id}G$, the associated Lie algebra.

About a century ago, Killing and Cartan classified complex simple
Lie algebras. Thanks to 
Coxeter and Dynkin, the classification can be expressed
pictorially. (See \cite{hawkins} for a wonderful account of their work
and the history surrounding it.) 

\begin{center}
\setlength{\unitlength}{3mm}
\begin{picture}(40,12)(2,-4)
\multiput(10,5)(2,0){6}{$\circ$}
\multiput(10.5,5.4)(2,0){5}{\line(1,0){1.6}} 
\put(14,3){$A_n$}

\multiput(25,5)(2,0){5}{$\circ$}
\multiput(25.5,5.4)(2,0){4}{\line(1,0){1.6}} 
\put(33.4,5.5){\line(1,1){1.2}} 
\put(33.45,5){\line(1,-1){1.2}}
\put(34.6,6.5){$\circ$}  
\put(34.6,3.4){$\circ$}  
\put(29,3){$D_n$}

\multiput(0,0)(2,0){5}{$\circ$}
\multiput(0.5,.4)(2,0){4}{\line(1,0){1.6}} 
\put(4,-2){$\circ$}
\put(4.3,-1.45){\line(0,1){1.5}}
\put(6,-2){$E_6$}

\multiput(15,0)(2,0){6}{$\circ$}
\multiput(15.5,.4)(2,0){5}{\line(1,0){1.6}} 
\put(19,-2){$\circ$}
\put(19.4,-1.45){\line(0,1){1.5}}
\put(21,-2){$E_7$}

\multiput(30,0)(2,0){7}{$\circ$}
\multiput(30.5,.4)(2,0){6}{\line(1,0){1.6}} 
\put(34,-2){$\circ$}
\put(34.3,-1.45){\line(0,1){1.5}}
\put(36,-2){$E_8$}

\end{picture} 
\end{center}

If a diagram has symmetry, we are allowed to fold it along the
symmetry and place an arrow pointing away from the hinge to get a new one:

\begin{center}
\setlength{\unitlength}{3mm}
\begin{picture}(40,10)(2,-2)
\multiput(10,5)(2,0){6}{$\circ$}
\multiput(10.5,5.4)(2,0){4}{\line(1,0){1.6}} 
\multiput(18.4,5.1)(0,.4){2}{\line(1,0){1.7}} 
\put(18.7,5){$>$}
\put(14,3){$B_n$}

\multiput(25,5)(2,0){6}{$\circ$}
\multiput(25.5,5.4)(2,0){4}{\line(1,0){1.6}} 
\multiput(33.4,5.1)(0,.4){2}{\line(1,0){1.7}} 
\put(33.8,5){$<$}
\put(29,3){$C_n$}

\multiput(10,0)(2,0){4}{$\circ$}
\multiput(10.5,.4)(4,0){2}{\line(1,0){1.6}} 
\multiput(12.4,.15)(0,.4){2}{\line(1,0){1.7}} 
\put(12.7,0){$>$}
\put(14,-2){$F_4$}

\multiput(28.5,0)(2,0){2}{$\circ$}
\multiput(28.9,.1)(0,.2){3}{\line(1,0){1.7}} 
\put(29.3,0){$>$}
\put(29,-2){$G_2$}

\end{picture} 
\end{center}
For example, $C_n$ is the fold of $A_{2n-1}$.
Given a semi-simple Lie algebra $\fg$, we let $D(\fg )$ denote
its Dynkin diagram.
\bigskip

Returning to geometry, homogeneous varieties correspond to marked 
diagrams

\begin{center}
\setlength{\unitlength}{3mm}
\begin{picture}(40,8)(2,1)
\multiput(10,5)(2,0){6}{$\circ$}
\multiput(10.5,5.4)(2,0){5}{\line(1,0){1.6}} 
\put(12,3){$G(3,7)$}
\put(14,5){$\bullet$}

\multiput(25,5)(2,0){5}{$\circ$}
\multiput(25.5,5.4)(2,0){4}{\line(1,0){1.6}} 
\put(33.4,5.5){\line(1,1){1.2}} 
\put(33.45,5){\line(1,-1){1.2}}
\put(34.6,6.5){$\circ$}  
\put(34.6,3.4){$\circ$}  
\put(27,3){$G_Q(3,14)$}
\put(29,5){$\bullet$}
\end{picture}
\end{center}

Note in particular that $\BP^n$ has marked diagram 

\begin{center}
\setlength{\unitlength}{3mm}
\begin{picture}(16,5)(3,2)
\multiput(6,5)(2,0){6}{$\circ$}
\multiput(6.5,5.4)(2,0){5}{\line(1,0){1.6}} 
\put(10,3){$\PP^n$}
\put(6,5){$\bullet$}
\end{picture}
\end{center}

 If $S$ is a subset of the simple roots of $G$,
 we use $P_S$ to designate the parabolic subgroup of $G$ corresponding to
the simple roots in the complement of $S$.  The Dynkin diagram for $G/P_S$ is the
diagram for $G$ with roots in $S$ marked.  Also, if $S$ consists of a single root
$\{\alpha\}$ and $\alpha=\alpha_i$,  we write $P_S=P_\alpha=P_i$.

\medskip

Unless otherwise specified, we   take $\ppp V$ to be 
ambient space for the the minimal
homogeneous embedding of $G/P$. For example, if $P=P_{\a}$ is
maximal, then $V=V_{\o}$ where $\o$ is the weight Killing-dual to
the coroot of $\a$.

\bigskip

Back to our question of how to study homogeneous varieties:

\medskip

\noindent{\bf Idea}: Study $X\subset\ppp V$ {\it via} $Y=Y_{X,x}\subset\ppp T_xX$.

\medskip

For this idea to be a good one, $Y$ should
be a simpler space than $X$  and it should
be possible to determine   $Y$ systematically.   
Although $Y$ is indeed cut out by quadratic equations, and has dimension
strictly less than that of $X$,  $X$ is  homogeneous,
and $Y$ need not be. 

The following theorem originates with work of Tits \cite{tits}, amplified
by Cohen and Cooperstein \cite{coco}:

\begin{theorem} \cite{LM0} Let $X=G/P_{\a}\subset\ppp V$ be
a rational homogeneous variety such that $\a$ is not
short (i.e., no arrow in $D(\fg)$ points towards $\a$).

Then $Y\subset\ppp T_xX$ is homogeneous, in fact  generalized minuscule
(as defined above).

Moreover, $Y$ can be determined pictorially:
remove the node corresponding to $\a$ from $D(\fg)$ and mark
the nodes that were adjacent to $\a$. One obtains a
semi-simple Lie algebra $\fh$ with marked diagram. The resulting
homogeneous space $H/Q$ is $Y$.
\end{theorem}

\begin{example} $ $

\begin{center}
\setlength{\unitlength}{3mm}
\begin{picture}(20,7)(13,.5)
\multiput(8,5)(2,0){5}{$\circ$}
\multiput(8.5,5.4)(2,0){4}{\line(1,0){1.6}} 
\put(16.4,5.5){\line(1,1){1.2}} 
\put(16.4,5){\line(1,-1){1.2}}
\put(17.6,6.5){$\circ$}  
\put(17.6,3.5){$\circ$}  
\put(8,2){$X=G_Q(4,12)$}
\put(14,5){$\bullet$}

\put(21,5){$\lra$}

\multiput(27,5)(2,0){2}{$\circ$}
\put(31,5){$\bullet$}
\multiput(27.5,5.4)(2,0){2}{\line(1,0){1.6}}
\put(34,5){$\bullet$}
\put(34.4,5.5){\line(1,1){1.2}} 
\put(34.4,5){\line(1,-1){1.2}}
\put(35.6,6.5){$\circ$}  
\put(35.6,3.5){$\circ$}  
\put(26,2){$Y=Seg(\pp 3\times G(2,4))$}
\end{picture}
\end{center}
\end{example}

The $H$-modules $\langle Y\rangle$ are
studied in \cite{kactheta}
where they are  called {\it type-I} $\theta$-representations.

The embedding of $Y$ is minimal iff the diagram is simply laced. In
the case of a double edge, one takes the quadratic Veronese
embedding, for a triple edge, one takes the cubic Veronese
(see the   algorithms below).

\smallskip

\begin{remark} 
In \cite{LM0} we explicitly determine $Y$ in the case of short roots and
in fact arbitrary parabolics. We also give geometric
models. In the case of $X=G/P_{\a}$ with
$\a$ short, $Y$ is the union of exactly two $G$-orbits, and the closed orbit has
codimension one in $Y$.\end{remark}

If $Y$ is still complicated, one can continue, studying the
asymptotic directions of $Y$ at a point. Eventually one gets (Segre products of re-embeddings of)
 $\BC\pp 1$, the one homogeneous space we all pretend to understand.
 
 We   describe our discovery that one can reverse this
 infinitesimalization
 procedure below. First we need to review some elementary algebraic
 geometry.

\subsubsection{Embeddings of projective
space}

Recall the Veronese embeddings of projective space:
$$
\begin{aligned}
v_d: \ppp V&\ra \ppp (S^dV) \\
[w]&\mapsto [w^d].\end{aligned}
$$
Dually, let $P_0\hd P_N$ be a basis of $S^dV^*$. The map  is
$$
[w]\mapsto [P_0(w)\hd P_N(w)].
$$

A remarkable fact is that all   maps between projective
spaces $\ppp V\ra \ppp W$ are obtained by projecting a Veronese
re-embedding.

If we project to $\ppp (S^dV/L)$ where $L\upperp\subset S^dV^*$ has
basis $Q_0\hd Q_m$, the map is
$$
\begin{aligned}
f: \ppp V&\ra \ppp (S^dV/L).\\
[w]&\mapsto [Q_0(w)\hd Q_m(w)]
\end{aligned}
$$

The image of $f$ is smooth and isomorphic to $\ppp V$ iff no secant
line of $v_d(\ppp V)$ intersects $L$.

$$\epsfbox{vog.6}$$

In the following algorithm we will actually be interested in images
that get squashed in the projection so that they are no longer
isomorphic to $\pp n$ (but are still smooth).

\subsection{First algorithm}

We are about to describe an algorithm, which you might like to think
of as a game. Starting with $\pp 1$
as initial input we  build some new
algebraic varieties subject to certain rules. The game has
rounds, and in each new round, we are allowed to use the outputs
from previous rounds as new inputs.
\medskip

We fix some notation.

Let $X\subset\ppp V$,
$Y\subset \ppp W$, be  varieties.

 $Seg(X\times Y)\subset\ppp (V\ot W)$ is
the Segre product, $([x],[y])\mapsto [x\ot y]$. Sometimes
we  just write   $X\times Y$ with the Segre product being
understood.

 $\s_k(X)\subset\ppp V$ is the union of all secant $\pp{k-1}$'s to 
$X$. We let $\s (X)=\s_2(X)$.

$\cT (X)\subset G(2,V)\subset\ppp(\La 2\ppp V)$ denotes the union of
all tangent lines to $V$. (Recall $G(k,V)$ is also the space of
$\pp{k-1}$'s in $\ppp V$.)

$\hat X\subset V$ is the cone over $X$ and
$\langle X\rangle\subseteq V$ denotes the linear span of $\hat X$.

Given $Y\subset \BP W$, we  let $I_k(Y)\subset S^kW^*$ denote the component
of the ideal in degree $k$. In the flowchart below we slightly abuse
notation by letting $I_k(Y)$  denote a set of generators of $I_k(Y)$.
Similarly, we let $T^*$ be shorthand for a basis of $T^*$.


\bigskip\centerline{{\bf\Large First Algorithm}}
\begin{center}
\setlength{\unitlength}{3mm}
\begin{picture}(30,48)(0,-5)
\put(-7.5,37.5){\line(1,0){7.8}}
\put(-7.5,41){\line(1,0){7.8}}
\put(-7.5,37.5){\line(0,1){3.5}}
\put(0.3,37.5){\line(0,1){3.5}}
\put(-7,39.5){Initial input}
\put(-6,38){$\PP^1\subset\PP^1$}
\put(1,38.5){\vector(1,0){3}}
\put(6,38){RESERVOIR}
\put(5.5,37.5){\line(1,0){8.5}}
\put(5.5,39.5){\line(1,0){8.5}}
\put(5.5,37.5){\line(0,1){2}}
\put(14,37.5){\line(0,1){2}}
\put(9.5,36.5){\vector(0,-1){3}}
\put(-2,31){Let $Y=Seg(v_{d_1}(X_1)\times\cdots\times v_{d_r}(X_r))
\subset\PP^{n-1}=\PP T$}
\put(0,29){with $X_j\subset\PP^{N_j}$ in the RESERVOIR}
\put(-2.5,28.5){\line(1,0){30}}
\put(-2.5,32.3){\line(1,0){30}}
\put(-2.5,28.5){\line(0,1){3.8}}
\put(27.5,28.5){\line(0,1){3.8}}
\put(9.5,27.5){\vector(0,-1){3}}
\put(3.5,23.5){\line(1,0){11.5}}
\put(4,22){Is $\langle\cT (Y)\rangle = \La 2T$ ?}
\put(3.5,21.5){\line(1,0){11.5}}
\put(15,21.5){\line(0,1){2}}
\put(3.5,21.5){\line(0,1){2}}
\put(8.5,21){\vector(-1,-1){6}}
\put(2.5,18){NO}

\put(-2.8,8.2){\line(1,0){3}}
\put(-2.8,8.2){\line(-1,3){2}}
\put(.2,8.2){\line(1,3){2}}
\put(-4.8,14.2){\line(1,0){7}}
\put(-3.5,12){TRASH}

\put(10,21){\vector(1,-1){6}}
\put(14.5,18){YES}
\put(9.5,14.5){\line(1,0){31}}
\put(10,13){Let $d$ be the smallest integer such that}
\put(20,11){$\s_d(Y)=\PP T = \PP^{n-1}$}
\put(10,9){Linearly embed $\PP T\subset\PP^n$ as the hyperplane $\{l=0\}$.}
\put(10,7.5){Consider the rational map}
\put(12,5){$\varphi : \PP^n\lra\PP^N\subset\PP(S^d\CC^{n+1})$}
\put(24.5,3.3){$||$}
\put(17,1.4){$\PP\langle l^d, l^{d-1}T^*, l^{d-2}I_2(Y), 
\ldots , I_d(\s_{d-1}(Y))\rangle .$}
\put(10,-1){Let $X=\overline{\varphi(\PP ^n)}\subset\PP^N$. 
Call $X$ an OUTPUT.}
\put(34,15.5){\line(0,1){23}}
\put(34,38.5){\vector(-1,0){19}}
\put(9.5,-2){\line(1,0){31}}
\put(9.5,-2){\line(0,1){16.5}}
\put(40.5,-2){\line(0,1){16.5}}
\end{picture}
\end{center}

For the first run through the algorithm, the admissible varieties and their
outputs are as follows
$$\begin{array}{rclcrcl}
Y & \subseteq & \pp{n-1}& & X^n & \subseteq & \pp N\\ \pp 1 & \subseteq &
\pp 1 & &\pp 2 & \subseteq & \pp 2\\ \pp 1\times\pp 1 & \subset & \pp 3 &
&\QQ^4 & \subset & \pp 5\\
 v_2(\pp 1) & \subset & \pp 2 & & \QQ^3 & \subset
& \pp 4. \\
\end{array} $$

Here and below,
$\QQ^m\subset\pp{m+1}$ denotes the smooth quadric hypersurface. 

For the second round,
$$\begin{array}{rclcrcl}
Y & \subseteq & \pp{n-1} & & X^n & \subseteq & \pp N\\ \pp 2 & \subseteq &
\pp 2 & &\pp 3 & \subseteq & \pp 3\\
v_2(\pp 2) & \subset & \pp 5 & & G_{\omega}(3,6) & \subset & \pp {11}\\ \pp
1\times\pp 2 & \subset & \pp 5 & &G(2,5) & \subset & \pp 9\\ \pp 2\times\pp
2 & \subset & \pp 8 & & G(3,6) & \subset & \pp{19}\\ \QQ^4 & \subset & \pp
5 & & \QQ^6 & \subset & \pp 8\\
\QQ^3 & \subset & \pp 4 & & \QQ^5 & \subset & \pp 7 .\end{array} $$
Here $G(k,l)$ denotes the Grassmannian of $k$-planes in $\BC^l$ and
$G_{\omega}(k,2k)$ denotes the
Grassmannian of Lagrangian $k$-planes for a given symplectic form.

\begin{question}What comes out?  Since the algorithm goes on forever, is it
even possible to answer this question?
\end{question}

\begin{proposition} \cite{LM1}
The algorithm is effective. A priori,    $r, d_j\leq 2$ and
the algorithm stabilizes after six rounds.
\end{proposition}

So at least our question is reasonable.
Now for the answer:

\begin{theorem} \cite{LM1} OUTPUTS = MINUSCULE VARIETIES.\end{theorem}

\begin{coro} \cite{LM1} A  new proof of the classification of the
compact Hermitian symmetric spaces without any reference to Lie groups.
\end{coro}

We have the following stable round:
$$\begin{array}{rclcrcl}
Y & \subseteq & \pp{n-1} & & X^n & \subseteq & \pp N\\ \pp {n-1} &
\subseteq & \pp {n-1} & &\pp n & \subseteq & \pp n\\ v_2(\pp {m-1}) &
\subset & \pp {\binom {m+1}2-1 } & & G_{\omega}(m, 2m) & \subset & \pp
{C_{m+1}-1}\\ \pp {k-1}\times\pp {l-1} & \subset & \pp{kl-1} & & G(k,k+l) &
\subset & \pp {\binom{k+l}{k}-1} \\
\QQ^{2m-2} & \subset & \pp {2m-1} & & \QQ^{2m} & \subset & \pp{2m+1}\\
\QQ^{2m-1}& \subset & \pp {2m} & &\QQ^{2m+1} & \subset & \pp{2m+2}\\ G(2,m)
& \subset & \pp{\binom m2 -1} & & \BS_m & \subset & \pp{2^{m-1}-1}.
\end{array} $$

Here $C_{m+1}=\frac{1}{m+2}\binom{2m+2}{m+1}$ is the $(m+1)$-st Catalan
number. The spinor variety $\BS_m$ of $D_m$ consists of one family of
maximal isotropic subspaces of $\CC^{2m}$ endowed with a nondegenerate quadratic form and embedded in the projectivization of one of
the two half-spin representations.

 The most interesting (but terminal) path is:
$$\begin{array}{rcl}
Seg(\pp 1\times\pp 2) & = & A_1/P_1\times A_2/P_2 \\  & \downarrow & \\
G(2,5) &= & A_4/P_2 \\ & \downarrow & \\ \BS_5 & = & D_5/P_5 \\ & \downarrow & \\
\BO\BP^2_\BC & = & E_6/P_6 \\ & \downarrow & \\  G_w(\BO^3,\BO^6) & = & E_7/P_7.
\end{array}$$

Here  $\OO\PP^2_{\CC}=E_{6(-14)}$ is the complexification of
  $\OO\PP^2=F_4/{\rm Spin}_9$,  
   the celebrated Cayley plane discovered by Chevalley. 
As a topological space, $\OO\PP^2$ is built out of
  three cells of dimension $0, 8$ and $16$.    
It is $F_{4(-20)}$ in the notation of Tits (see e.g. \cite{loos}).  
The notation  $ G_w(\BO^3,\BO^6)$ is discussed in \S 4.6.
 
 \medskip
 
 \subsection{Second algorithm}
 
 Now that was fun, but it was a shame to throw away some of our
 favorite varieties like $v_3(\pp 1)$ and $Seg(\pp 1\times\pp 1\times\pp 1)$.
Let's   revise the algorithm slightly and have just one
consolation round. Given 
a smooth variety $Y\subset \BP V$, we let $\tau (Y)\subset \BP V$
denote the union of all points on all embedded tangent lines to $Y$. 

\bigskip
\
\bigskip


\centerline{{\bf\Large Consolation prize (one round)}}
\begin{center}
\setlength{\unitlength}{3mm}
\begin{picture}(30,43)(5,-1)
\put(1,38){MINUSCULE VARIETIES}
\put(.5,37.5){\line(1,0){16.5}}
\put(.5,39.5){\line(1,0){16.5}}
\put(.5,37.5){\line(0,1){2}}
\put(17,37.5){\line(0,1){2}}
\put(9.5,36.5){\vector(0,-1){3}}
\put(-2,31){Let $Y=Seg(v_{d_1}(X_1)\times\cdots\times v_{d_r}(X_r))
\subset\PP^{n-2}=\PP T_1$}
\put(0,29){with $X_j\subset\PP^{N_j}$ minuscule}
\put(-2.5,28.5){\line(1,0){30.3}}
\put(-2.5,32.3){\line(1,0){30.3}}
\put(-2.5,28.5){\line(0,1){3.8}}
\put(27.8,28.5){\line(0,1){3.8}}
\put(9.5,27.5){\vector(0,-1){3}}
\put(2.5,23.5){\line(1,0){13.5}}
\put(3,22){Is codim $\langle\cT (Y)\rangle = 1$  ?}
\put(2.5,21.5){\line(1,0){13.5}}
\put(16,21.5){\line(0,1){2}}
\put(2.5,21.5){\line(0,1){2}}
\put(8.5,21){\vector(-1,-1){6}}
\put(2.5,18){NO}

\put(-2.8,8.2){\line(1,0){3}}
\put(-2.8,8.2){\line(-1,3){2}}
\put(.2,8.2){\line(1,3){2}}
\put(-4.8,14.2){\line(1,0){7}}
\put(-3.5,12){TRASH}

\put(10,21){\vector(1,-1){6}}
\put(14.5,18){YES}
\put(9.5,14.5){\line(1,0){37}}
\put(10,13){Linearly embed $\PP T_1\subset_{\{l=0\}}\PP^{n-1}\subset_{\{m=0\}}\PP^n$.}
\put(10,11.5){Consider the rational map}
\put(12,9){$\varphi : \PP^n\lra\PP^N\subset\PP(S^4\CC^{n+1})$}
\put(24.5,7.3){$||$}
\put(10,5.4){$\PP\langle l^4,l^3m, l^3T_1^*, l^2I_2(Y), 
l^2mT_1^*-lI_3(\tau(Y)_{sing}),
l^2m^2-I_4(\tau(Y))\rangle .$}
\put(10,3.5){Let $X=\overline{\varphi(\PP ^n)}\subset\PP^N$. 
Call $X$ an OUTPUT.}
\put(9.5,3){\line(1,0){37}}
\put(9.5,3){\line(0,1){11.5}}
\put(46.5,3){\line(0,1){11.5}}
\put(25,2){\vector(0,-1){3}}
\put(22,-3){OUTPUTS}
\put(21.5,-1.5){\line(1,0){7.5}}
\put(21.5,-3.5){\line(1,0){7.5}}
\put(21.5,-3.5){\line(0,1){2}}
\put(29,-3.5){\line(0,1){2}}
\end{picture}
\end{center}

\bigskip

\bigskip

\bigskip

This algorithm is also effective. We show   
{\it a priori} that $r,d_j\leq 3$.

So, what do we get???

\begin{theorem} \cite{LM1} OUTPUTS = Fundamental adjoint varieties $X\subset\ppp \fg$.
\end{theorem}

The fundamental adjoint varieties are the closed orbits in the adjoint
representation when it is fundamental. The Lie algebras whose
adjoint varieties are not fundamental are the pathological 
$A_n$ and its fold (when foldable) $C_m$. Our algorithm again provides a beautiful
and easy proof of the classification of fundamental adjoint varieties.
We   also  account for these
two pathological cases, but it involves some less elegant work.
We do obtain:

\begin{coro} \cite{LM1}
A new proof of the Cartan-Killing classification
of complex simple Lie algebras.\end{coro}

\pagebreak
The examples for the adjoint algorithm are as follows:

$$\begin{array}{rclcc}
Y & \subset & \PP^{n-2}& &\fg\\
v_3(\pp 1)& \subset & \pp 3 & & \fg_2\\
\pp 1\times \QQ^{m-4} & \subset & \pp {2m-5}& & \fso_m\\ 
G_{\omega}(3,6) & \subset & \pp
{13}& & \ff_4\\ G(3,6)& \subset & \pp {19}& & \fe_6\\
\BS_6& \subset & \pp {31}& & \fe_7\\
G_w(\BO^3, \BO^6)& \subset & \pp {55}& & \fe_8 .\end{array}$$
The two exceptional (i.e., non-fundamental) cases are $$\begin{array}{rclcc}
\pp {k-3}\sqcup\pp {k-3} & \subset & \pp{2k-3} & & \fsl_k\\ \emptyset &
\subset & \pp {2m-1} & & \fsp_{2m} . \end{array} $$

Note that the algorithm works in
these two cases, we just didn't have the varieties in
our reservoir. (In the case of $\fsp_{2m}$, it may have been there
and just difficult to see.)

\subsection{Outline of the proofs}
The proofs have three ingredients: differential invariants, local Lie
algebras, and relating Casimirs to geometry.

\subsubsection { Differential invariants} Given a variety $X\subset\ppp V$ and $x\in X$,
we can recover $X$ from its Taylor series at $x$. The projective differential
invariants comprise a series of tensors encoding the geometric (i.e., invariant
under $GL(V)$) information in the Taylor series. We prove {\it a priori}
facts about the differential invariants of any putative minuscule
variety or fundamental adjoint variety.

\subsubsection { Local Lie algebras} 
A local Lie algebra is a graded vector space
$$
\fg_{-1}\op\fg_0\op\fg_1
$$
equipped with a bracket for which
$\fg_0$ is a Lie algebra.  The bracket must respect the grading but the Jacobi
identity need not hold. 
  If the Jacobi identity fails, one can construct a (unique) $\ZZ$-graded
Lie algebra from the local Lie algebra. The traditional way to
do this (see, e.g., \cite{kac}) is to take the free algebra
generated by the brackets and then mod out by the relations. Note
that if one does that, one has no idea how many factors (if any
at all) one will be adding on to obtain the final result.

The data $Y=H/Q\subseteq\ppp T_1$
furnishes (up to scale) a local Lie algebra  
with $\fg_0=\fh\op\BC$, $\fg_1=T_1$, $\fg_{-1}=T_1^*$.
Here the action of $\BC$ is as a scalar times the identity and we
do not initially specify the scalar.

Since this does give rise to a unique $\ZZ$-graded Lie algebra,
(in particular, a Lie algebra equipped with a representation $V$ supported
on one fundamental weight), we can study the resulting homogeneous
variety $  X'=G/P\subset\ppp V$ , and calculate its differential
invariants.

Note that if $   X'$ is minuscule, $\fg =\fg_{-1}\op\fg_0\op\fg_1$
and if $  X'$ is adjoint,  $\fg=\fg_{-2}\op\fg_{-1}\op\fg_0\op\fg_1
\op\fg_2$ with $\fg_{\pm 2}=\BC$.

\subsubsection {Compare $  X'$ with the constructed variety $X$}
We compare the differential invariants of $X$ and $  X'$ and show that
they are the same.
The key point in the
minuscule
case  is that a local Lie algebra is already a Lie algebra
iff the Jacobi identity holds. We show this is the case
iff $\La 2\fg_1$ is an eigenspace for the Casimir operator of $\fh$.
We then connect this to the geometry by showing that 
$\langle \cT (Y)\rangle$ is a Casimir eigenspace!

For the adjoint case the idea is similar: one shows that one can get
away with a one-dimensional correction.

 \pagebreak

\subsection{Applications, generalizations and related work}

\subsubsection{More general algorithms} The adjoint algorithm
(with a different test for admission) also constructs all
homogeneous varieties and representations corresponding
to a five step grading. M. Dillon is currently
  formulating more general constructions for
  gradings including extensions to affine Lie algebras.
(The algorithms have no hope of being effective for graded Lie
algebras with exponential growth.)
An amusing exercise is to construct the grading of $\fe_8\up 1$
associated to $\a_5$. Here $\fh= \fa_4+\fa_4$, and 
$$\begin{array}{rccclcl} Y=Y_1 & = & G(2,5)\times
\pp 4=Seg(G(2,V)\times \ppp W) & \subset  & \ppp (\La 2 V\ot W) & = & \ppp\fg_1, \\
Y_2 & = & Seg(\ppp V^*\times G(2,W)) & \subset & 
\ppp (V^*\ot \La 2 W) & = & \ppp\fg_2, \\ 
Y_3 & = & Seg(\ppp V \times G(2,W^*)) & \subset & 
\ppp (V \ot \La 2 W^*)  & = & \ppp\fg_3, \\ 
Y_4 & = & Seg(G(2,V)\times \ppp W^*) & \subset & \ppp (\La 2 V\ot W^*) & = & \ppp\fg_4,
\end{array}$$
and then the cycle repeats. Note that all the varieties are isomorphic
as varieties and they appear in all possible combinations in terms
of the $\fh$-action.

 In this context, a result of Kostant \cite{koscubic},
 (theorem 1.50) is interesting.
He determines, given a reductive Lie algebra $\fr$ and an $\fr$-module
$V$, when $\fr + V$ can be given the structure of a Lie algebra
compatible with the $\fr$-actions (and in how many different ways).

\subsubsection{Constructions in an algebraic context}
Let $\cA$ be an algebra with  unit and involution $a\mapsto \overline a$.
B. Allison \cite{al} defines $\cA$ to be {\it structurable} if the following  
graded vector space is actually a graded Lie algebra. Let
$$\fg (\cA) = \fg_{-2}\op\fg_{-1}\op\fg_0\op\fg_1\op\fg_2
$$
with $\fg_{\pm 1}\simeq \cA$, $\fg_{\pm 2}\simeq \cS (\cA):= \{
a\in\cA \mid \overline a=-a\}$ and $\fg_0$ the set of derivations
of $\cA$ generated by linear maps $V_{a,b} : \cA\ra \cA$, where
$V_{a,b}(c)=(a\overline b)c+(c\overline b)a - (c\overline a)b$.

$\fg (\cA)$ always has a natural bracket, and is a Lie algebra
iff the Jacobi identity holds. In analogy with our situation, there
is essentially one identity to check, which, in Allison's notation
becomes
$$
[V_{a,1},V_{b,c}]=V_{V_{a,1}(b),c}- V_{b,V_{\overline a,1}(c)}.
$$
 Allison shows there is a one to one correspondence between
  simple structurable algebras and certain symmetric five step gradings of
 simple Lie algebras. His constructions
 work over all  fields and one motivation for his constructions
 was to determine simple Lie algebras over arbitrary fields.
 Our more general $5$-step constructions correspond to algebraic structures
termed
 {\it Kantor triple systems}.

{\it Jordan triple systems} form a special class of Kantor triple systems,
those with trivial involution $\overline a=a$ and thus   $\cS (\cA)=0$,
and they give rise to the minuscule gradings our 
minuscule algorithm  produces.

\subsubsection{Rigidity}
In contrast to  the rigidity theorems in \S 2.1.1  we have:

\begin{theorem} There exist \lq\lq fake\rq\rq\ adjoint varieties.
That is, for each adjoint variety $X\subset\ppp \fg$, there
exists $Z\subset\pp N$, not isomorphic to $X$,   
$U\subset Z$ a Zariski open subset, and a holomorphic
map $\phi: U\ra X$ such that asymptotic directions are preserved
 (in fact the
entire projective second fundamental form is preserved,
see \cite{LM1}).

Moreover, the same is true for all non-minuscule homogeneous varieties.
\end{theorem}

This local flexibility fails globally for adjoint varieties.
Consider the following result of Hong \cite{hong}:

\begin{theorem} A Fano manifold with a geometric structure modeled after
a fundamental adjoint variety $Z$ is
biholomorphic to $Z$ and the geometric structure is locally isomorphic
to the standard geometric structure
on $Z$.
\end{theorem}

A geometric structure on a variety $X$ modeled after
a fundamental adjoint variety $Z$ may be understood as
follows. Tor all
$x\in X$, one has a subvariety $Y_x\subset \BP T_xX$ isomorphic
to the asymptotic directions $Y_z\subset\BP T_zZ$. Note that here
the subvariety $Y_x$ is not required to play any particular role,
but it must be present at every point so
in  particular  it determines a 
family of distributions
on $X$.

Hong's theorem is a variant of an earlier rigidity result of Hwang
and Mok \cite{hmrigid}, where a geometric structure may
be understood in the analogous way.

\begin{theorem} A Fano manifold with a geometric structure modeled after
a compact irreducible Hermitian symmetric space $S$ of rank $\geq 2$ is
biholomorphic to $S$ and the geometric structure is locally isomorphic
to the standard geometric structure
on $S$.
\end{theorem}

The work of Hwang and Mok (also see \cite{hwangmokpapers}) relies
on studying an intrinsic analog of the set of asymptotic directions.
Namely a Fano manifold $X$ is uniruled by rational curves, and, fixing a
reference line bundle, one has a subvariety $Y_x\subset \BP T_xX$
of tangent directions to minimal degree rational curves.  
They also study deformation rigidity of homogeneous varieties with Picard number one. 
It is worth remarking that in 
the case of Hermitian symmetric spaces \cite{hmhss}, a key fact used in their
proof is that $\cT (Y_x)\subset\BP T_xX$ is linearly nondegenerate, i.e.,
the same condition we use in the minuscule algorithm.

\medskip

\subsubsection{Normal forms for singularities}
In \cite{Ar}, Arnold classified the simple singularities.
Like many interesting things in life they are in correspondence
with 
Dynkin diagrams, in fact just the simply laced ones. He also gave normal
forms in a minimal set of variables (two). It was also known (but
evidently unpublished) that the simple singularities can be realized
as degree three hypersurface singularities if one allows the number
of variables to grow with the Milnor number. Holweck \cite{holweck}
has found a nice realization of these hypersurface singularities using
a theorem of Knop \cite{knop} and the construction in \S 2.3.

\subsection{Why secant and tangent lines?}

The idea that secant and tangent lines should so strongly control
the geometry of homogeneous varieties was inspired by 
Zak's theorem on Severi varieties.

\begin{theorem} [Zak's theorems on linear normality and Severi varieties]
Let $X^n\subset\pp\na$ be a smooth variety,
not contained in a hyperplane and such that $\s (X)\neq\pp\na$
Then

i. $a\geq \frac n2+2$.

ii. If $a=\frac n2+2$ then $X$ is one of $v_2(\pp 2)\subset\pp 5$,
$Seg(\pp 2\times\pp 2)\subset\pp 8$, $G(2,6)\subset\pp{14}$,
${\BO\BP^2_\BC}\subset\pp{26}$.
\end{theorem}

The four critical varieties are called the {\it Severi} varieties
after Severi who proved the $n=2$ case of the theorem.

It is not known if there is a bound on the secant defect
$\delta_{\s}(X):= 2n+1-\tdim\,\s (X)$ for
smooth subvarieties of projective space with degenerate
secant varieties. On the other hand,
  Zak established an upper bound
  on the codimension of a smooth variety
  of a given secant defect. He then went on
to classify the varieties achieving this bound, which he calls the
{\it Scorza varieties}.
They are all closed orbits $G/P\subset\ppp V$,
namely $v_2(\pp n)$, $Seg(\pp n\times \pp n)$, $G(2,n)$
and ${\BO\BP^2_\BC}$.

Using Zak's result, Ein and Shepherd-Baron proved these four
varieties also classify quadro-quadro Cremona transformations
\cite{ESB}.
There are numerous other characterization problems where the
answer is the Severi varieties.

Zak's proofs of his theorems rely on looking at 
{\it entry loci}. Namely let $X\subset\ppp V$ be a variety
and let $y\in \ppp V\backslash X$. Define the entry locus of $y$ to be
$$
\Sigma_y:=\overline{ \{ x\in X\mid \exists z\in X, y\in \pp 1_{xz}\} }.
$$
 Here $\pp 1_{xz}$ denotes the projective line spanned
by $x$ and $z$.  Zak shows that for a Severi variety, the
entry locus of a general point of $\s (X)$ is a quadric hypersurface
in $X$. He then goes on to show that a variety so uniruled by quadric
hypersurfaces and satisfying the dimension requirements must satisfy
further dimension restrictions and eventually must be one of the four
Severi varieties. More precisely, he shows that each Severi variety
is the image of a rational map of a projective space. Our algorithms
  generalize his construction of the Severi varieties.

 Recently Chaput \cite{chaput1} has shown that such
a uniruling by quadrics immediately implies homogeneity, which gives
a quicker proof of Zak's theorem. 

The rigidity results for the Severi
varieties  other than $v_2(\pp 2)$,
are an outgrowth of a different proof of Zak's theorem,
where one first shows that any putative Severi variety infinitesimally
looks like an actual Severi variety to second order at a general point,
and then uses the rigidity to finish the proof. ($v_2(\pp 2)$
requires special treatment as its tangential variety is nondegenerate.)
  See \cite{Lsec} and \cite{Lrigid}.

\medskip Many other important problems are related to the geometry of 
secant varieties. Consider the classical Waring problem: given a generic
homogeneous polynomial of degree
$d$ in $n+1$ variables, what is the minimal $k$ such that the polynomial may
be written as a sum of $k$ $d$-th powers?
Phrased geometrically, the problem
is to find the minimal $k$ such that $ \s_k(v_d(\pp n))=\ppp (S^d \BC^{n+1})$.
This problem was solved by Alexander and Hirschowitz  \cite{ah}.
Generalizations and variants are still open, see 
\cite{cach, iaka, rasc} for recent progress.

\subsection{Tits correspondences and applications}

In \S 2.7.1 we describe a construction, {\it Tits correspondences} \cite{tits},
to determine the homogeneous  (in the 
sense described below) unirulings of a
rational homogeneous variety $X\subset\BP V$.   We explain how to use
Tits correspondences to explicitly construct the homogeneous
vector bundles over $X$ in \S 2.7.2. In \S 2.7.3 we relate Tits 
correspondences to other orbit closures in $\BP V$ and
  systematize Kempf's method for desingularizing
orbit closures.  Tits correspondences led us to the decomposition
formulas described in \S 4.2 and we believe they will have more
applications in the future.

\subsubsection{Tits correspondences}

 Tits associated to any simple group a full set of geometries, encoded by the
parabolic subgroups and their relative positions. 
This culminated in the definition of 
{\em buildings}, which has since known  formidable developments.
The Tits correspondences we now describe come  from this perspective.  

Let $G$ be a simple Lie group,
let $S,S'$ be two subsets   of simple roots of
$G$. Consider the diagram 
$$\begin{array}{ccccc} & & G/P_{S\cup S'} & & \\ &
{\scriptstyle \pi} \swarrow & &
\searrow {\scriptstyle \pi '} & \\
X= G/P_S & & & & X'= G/P_{S'}
\end{array}$$
Let $x'\in X'$ and consider $Y:= \pi (\pi'\inv (x'))\subset X$.
We call $Y=Y_{x'}$ the {\it Tits transform} of $x'$. The variety $X$ is
covered by such varieties $Y$. Tits shows that $Y= H/Q$ where
${\mathcal D}(H)={\mathcal D}(G)\backslash (S\backslash S')$, and 
$Q\subset H$ is the parabolic subgroup corresponding to $S'\backslash S$.
 We   call such subvarieties $Y$ of $X$, 
{\it $G$-homogeneous subvarieties}. 
 
\subsubsection{Unirulings}
 
$  $

\medskip\noindent {\sl Unirulings by lines}.
Let $X=G/P_{\a}\subset\ppp V$ be a homogeneous variety. In Theorem 2.5 above
we determined the parameter space of lines through a point of $X$. Here
we determine the space of all lines on $X$.
Let $\BF_1(X)^G\subset G(2,V)\subset\ppp\La 2V$ denote the variety 
parametrizing the $G$-homogeneous $\pp 1$'s on $X$.
Then $\BF_1(X)^G=G/P_{S}$ where $S$ is the set
of simple roots adjacent to $\a$ in $D(\fg)$.
If $P$ is not maximal, there is one family of $G$-homogeneous
$\pp 1$'s for each marked node.
If $\a$ is not a short root, all lines are $G$-homogeneous.

\medskip\noindent {\sl Unirulings by $\pp k$'s}.
The diagram for $\pp k$ is $(\fa_k,\o_1)$.
Here there may be several $G$-homogeneous families of $\pp k$'s on
$X=G/P_{\a}$, each arising from a subdiagram of $D(\fg, \o)$ isomorphic
to a $D(\fa_k,\o_1)$. If $\a$ is not short, then all $\pp k$'s
on $X$ are $G$-homogeneous.

\begin{example} The largest linear space on $E_n/P_1$ is a $\pp{n-1}$,
via the chain terminating with $\a_n$, so $E_n/P_1$ is  
maximally uniruled
by $\pp{n-1}$'s. There is a second chain terminating with
$\a_2$, so $E_n/P_1$ is also maximally uniruled by $\pp 4$'s. (The 
unirulings by the $\pp 4$'s are maximal in the sense that none
of the $\pp 4$'s of the uniruling are contained in any  
$\pp{5}$ on $E_n/P_1$.)  The varieties parametrizing these two unirulings
of $E_n/P_1$ are
respectively $E_n/P_2$ and $E_n/P_5$.  

\medskip

\setlength{\unitlength}{3mm}
\begin{picture}(20,3.5)(-10,-1.8)
\multiput(0,0)(2,0){6}{$\bullet$}
\multiput(0.45,.3)(2,0){5}{\line(1,0){1.65}}
\put(8,0){$\bullet$}
\put(4,-2){$\circ$}
\put(4,-2){$\ast$}
\put(4.3,-1.4){\line(0,1){1.5}}
\put(-.5,1.1){$\a_1$}

\multiput(16,0)(2,0){3}{$\bullet$}
\multiput(22,0)(2,0){3}{$\circ$}
\multiput(16.45,.3)(2,0){5}{\line(1,0){1.65}}
\put(20,-2){$\bullet$}
\put(22,0){$\ast$}
\put(20.3,-1.4){\line(0,1){1.5}}
\put(15.5,1.1){$\a_1$}
\end{picture}
\end{example}

\medskip\noindent {\sl Unirulings by quadrics}. 
Here we look for subdiagrams isomorphic to the diagram of a quadric
hypersurface, i.e.  $D(\fso_n,\o_1)$. Note that
there are at most two
possible such subdiagrams
except for  the case
of $(\fg,V)=(\fso_8,V_{\o_2})$
where there are three isomorphic ones.

\begin{example} $(\fg, V)=(\fso_{2m},V_{\o_k})$ has two such unirulings,
by $\BQ^4$'s parametrized by $D_m/P_{\a_{k-2},\a_{k+2}}$ and by
$\BQ^{2(m-k-1)}$'s parametrized by $D_m/P_{\a_{k-1} }$.

\begin{center}
\setlength{\unitlength}{3mm}
\begin{picture}(20,7)(11,1.5)
\multiput(8,5)(2,0){5}{$\circ$}
\multiput(8.5,5.4)(2,0){4}{\line(1,0){1.6}} 
\put(16.4,5.5){\line(1,1){1.2}} 
\put(16.4,5){\line(1,-1){1.2}}
\put(17.6,6.5){$\circ$}  
\put(17.6,3.5){$\circ$}  
\put(12,5){$\bullet$}
\put(11.5,4){$\a_k$}
\multiput(8,5)(8,0){2}{$\ast$}

\multiput(23,5)(2,0){5}{$\circ$}
\multiput(23.5,5.4)(2,0){4}{\line(1,0){1.6}} 
\put(31.4,5.5){\line(1,1){1.2}} 
\put(31.4,5){\line(1,-1){1.2}}
\put(32.6,6.5){$\circ$}  
\put(32.6,3.5){$\circ$}  
\put(27,5){$\bullet$}
\put(26.5,4){$\a_k$}
\put(25,5){$\ast$}

\end{picture}
\end{center}

\end{example}

\begin{example}[The Severi varieties] The
$G$-homogeneous quadrics uniruling
the Severi varieties are parametrized by the same Severi variety
in the dual projective space.

\begin{center}
\setlength{\unitlength}{3mm}
\begin{picture}(40,6)(2,2)
\multiput(10,5)(2,0){5}{$\circ$}
\multiput(10.5,5.4)(2,0){4}{\line(1,0){1.6}} 
\put(12,5){$\bullet$}
\put(16,5){$\ast$}
\put(12,3){$G(2,6)$}

\multiput(25,5)(2,0){5}{$\circ$}
\multiput(25.5,5.4)(2,0){4}{\line(1,0){1.6}} 
\put(29,3){$\circ$}
\put(25,5){$\bullet$}
\put(33,5){$\ast$}
\put(29.3,3.55){\line(0,1){1.5}}
\put(31,3){$\OO\PP^2_{\CC}$}

\end{picture}
\end{center}
\end{example}

Unirulings by quadrics give rise to subvarieties of $\s (X)$
as follows.
The union of the spans of the quadrics produce a subvariety
of $\s (X)$ whose entry loci contain the quadrics (and in all examples
we know, are equal to the quadrics). It also appears that,
when $X$
is homogeneous, these
give the subvarieties of $\s (X)$ with maximal entry loci. As explained below, these orbit closures admit
uniform desingularizations by Kempf's method.

\subsubsection{Homogeneous vector bundles}
The $G$-homogeneous vector bundles over 
rational homogeneous varieties $Z=G/Q$ are defined by
$Q$-modules. (We change notation to reserve
$G/P$ for a different role below.)
 If $W$ is a $Q$-module, one obtains the vector
bundle $E_W:=G\times_Q W\ra G/Q$ where $(gp,w)\simeq (g,pw)$ for all $p\in Q$.
Note that $D(\ff)$ is $D(\fg)$ with the nodes
corresponding to $Q$ deleted. 

\smallskip

Some homogeneous vector bundles over
rational homogeneous varieties $Z=G/Q$ can be
understood in terms of Tits fibrations.
Let $X=G/P\subset \BP V$. For each $z\in Z$ we obtain a variety
$Y_z\subset X$ and thus a linear space $\langle Y_z\rangle \subset V$.
As we vary $z$, we obtain a vector bundle $E\ra Z$ whose fibers
are the $\langle Y_z\rangle$'s. In particular, $E$ is a subbundle
of the trivial bundle $V \ot\mathcal O_{Z}\ra Z$.

\subsubsection{Collapses \`a la Kempf}

Let $X=G/P\subset \ppp V$ be a rational homogeneous variety.
In \cite{kempf}, Kempf defined a method for desingularizing
 orbit closures $\overline\cO \subset\ppp V$ by finding a vector bundle
over a different homogeneous space $Z=G/Q$. He calls this
technique
the {\it collapsing of a vector bundle}. He gave some examples
that appear to be found    ad hoc. In
\cite{LMmagic} we gave new examples of collapsing occurring in
series.
Using the above discussion, we now give another
description of Kempf's method.

\subsubsection*{Recipe for Kempf's method}
Let $V$ be a $G$-module with closed orbit $X\subset \BP V$.
Let $Z=G/Q$ be a rational homogeneous variety of $G$, defining a Tits correspondence with $X$.

 Consider the $G$-variety, which   often appears to be an orbit closure
$$
\overline\cO =\cup_{z\in Z}\langle Y_{z}\rangle \subset \BP V .
$$
If $\overline\cO\neq \ppp V$, then $\overline\cO$ is a singular variety.
It admits a desingularization by   $\ppp E$, where $E\ra Z$ is
the vector bundle constructed above
 as the natural
map $\ppp E\ra Z$ is generically one to one.

Note that  the $G$-variety
$\tau (X)$ is of course desingularized by $\ppp TX$, which is
 a special case of the above with $Q=P$. 
\bigskip

\def\hd{,...,}
\def\trank{\text{rank}}
\def\tdim{\text{dim}}
\def\fa{{\mathfrak a}}
\def\pp#1{\mathbb P^{#1}}

\def\com{{\rm com}}
\def\trace{{\rm trace}}
\def\cI{{\mathcal I}}
\def\cE{{\mathcal E}}
\def\cA{{\mathcal A}}
\def\cF{{\mathcal F}}
\def\cG{{\mathcal G}}
\def\cD{{\mathcal D}}
\def\cJ{{\mathcal J}}
\def\cZ{{\mathcal Z}}
\def\cR{{\mathcal R}}
\def\cS{{\mathcal S}}
\def\cL{{\mathcal L}}
\def\cW{{\mathcal W}}
\def\cO{{\mathcal O}}
\def\cQ{{\mathcal Q}}
\def\CC{\mathbb C}
\def\RR{\mathbb R}
\def\HH{\mathbb H}
\def\AA{{\mathbb A}}
\def\BB{{\mathbb B}}
\def\OO{\mathbb O}
\def\SC{\tilde{\mathbb C}}
\def\SR{\tilde{\mathbb R}}
\def\SH{\tilde{\mathbb H}}
\def\SA{\tilde{\mathbb A}}
\def\SB{\tilde{\mathbb B}}
\def\SO{\tilde{\mathbb O}}
\def\LG{\mathbb {LG}}
\def\LF{{\mathbb {LF}}}
\def\ZZ{\mathbb Z}
\def\SS{\mathbb S}
\def\GG{\mathbb G}
\def\11{\mathbf 1}
\def\PP{\mathbb P}
\def\QQ{\mathbb Q}
\def\FF{\mathbb F}
\def\JA{{\mathcal J}_3(\AA)}
\def\JB{{\mathcal J}_3(\BB)}
\def\ZA{{\mathcal Z}_2(\AA)}
\def\fh{{\mathfrak h}}
\def\fs{{\mathfrak s}}
\def\cS{{\mathfrak S}}
\def\fb{{\mathfrak b}}
\def\fc{{\mathfrak c}}
\def\fd{{\mathfrak d}}
\def\fsl{{\mathfrak {sl}}}
\def\fsu{{\mathfrak {su}}}
\def\fsp{{\mathfrak {sp}}}
\def\fspin{{\mathfrak {spin}}}
\def\fso{{\mathfrak {so}}}
\def\fe{{\mathfrak e}}
\def\ff{{\mathfrak f}}
\def\fz{{\mathfrak z}}
\def\ffi{{\mathfrak i}}
\def\fg{{\mathfrak g}}
\def\fn{{\mathfrak n}}
\def\fp{{\mathfrak p}}
\def\fk{{\mathfrak k}}
\def\ft{{\mathfrak t}}
\def\fl{{\mathfrak l}}
\def\etimes{\hspace{-1mm}\times\hspace{-1mm}}
\def\l{\lambda}
\def\a{\alpha}
\def\ta{\tilde{\alpha}}
\def\o{\omega}
\def\oo{\Omega}
\def\O{\Omega}
\def\b{\beta}
\def\g{\gamma}
\def\s{\sigma}
\def\k{\kappa}
\def\d{\delta}
\def\th{\theta}
\def\m{\mu}
\def\up#1{{}^{({#1})}}
\def\e{\varepsilon}
\def\ot{{\mathord{\,\otimes }\,}}
\def\op{{\mathord{\,\oplus }\,}}
\def\otc{{\mathord{\otimes\cdots\otimes}\;}}
\def\pc{{\mathord{+\cdots +}}}
\def\lra{{\mathord{\;\longrightarrow\;}}}
\def\ra{{\mathord{\;\rightarrow\;}}}
\def\da{{\mathord{\downarrow}}}
\def\we{{\mathord{{\scriptstyle \wedge}}}}
\def\JA{{\mathcal J}_3(\AA)}
\def\JB{{\mathcal J}_3(\BB)}
\def\tr{{\rm trace}\;}
\def\dim{{\rm dim}\;}
\def\La{\Lambda}
 
\section{Triality and exceptional Lie algebras}

\smallskip
We begin with a review of division algebras via
the Cayley-Dickson process in \S 3.1 and 
explicitly describe the  derivations of $\BO$
in \S 3.2 for later use. In \S 3.3 we review
the triality principle for the octonions and its extension
to all structurable algebras. We pause for a detour in \S 3.4,
exploring the triality model for $\fso_{4,4}$ in detail, which
ends up being a $4$-ality. We review various constructions of the
magic square in \S 3.5. We present a new result
regarding the compatiblity of the Cayley-Dickson process
and inclusions of triality algebras in \S 3.6, which
at first glance appears to be trivial, but whose \lq\lq obvious\rq\rq\
generalization to structurable algebras is false.
In \S 3.7 we study automorphisms related to the algebras in the magic
square and the resulting magic squares of symmetric and trisymmetric
spaces. In \S 3.8 we use our result in \S 3.6 to show how some
towers of dual pairs naturally occur in series. In \S 3.9 we discuss
quaternionic symmetric spaces, their relation to adjoint varieties
and the conjecture of LeBrun and Salamon.

\subsection{Division algebras} 

A classical theorem of Hurwitz, published in 1898, asserts that there are only 
four normed division algebras
 over the field of real numbers : $\RR$ itself, $\CC$,
the quaternions $\HH$, and the octonions $\OO$.  A {\em division 
algebra} is   an algebra without zero divisor  and    a {\it
normed division algebra} 
  is a division algebra endowed with a norm 
such that the norm of a product is the product of the norms. 

\smallskip
Quaternions were discovered by Hamilton in 1843, and octonions shortly afterwards 
by Graves and Cayley. A nice way to define these algebras is to imitate the 
definition of complex numbers by a pair of real numbers. Let $\AA$ be 
a real algebra endowed with a conjugation $x\mapsto \overline{x} $,
 such that $x+\overline{x}$ and 
$x\overline{x}=\overline{x}x$ are always real numbers (more precisely scalar multiples of the unit 
element), the later being positive when $x$ is nonzero.
If $\AA$ is alternative, i.e., 
  any subalgebra generated by two elements is associative \cite{baez},
then $\AA$ is a normed algebra 
for the norm $\norm{x}^2=x\overline{x}$.
 
\smallskip
One can then define  a new algebra with conjugation  $\BB=\AA\op\AA$ by letting
$(x,y)(z,t)=(xz-t\overline{y},\overline{x}t+zy)$ 
with conjugation  $\overline{(x,y)}=(\overline{x},-y)$. This is   the 
{\em Cayley-Dickson} process. The new algebra $\BB$ will
  be alternative,   hence a normed algebra, iff $\AA$ is 
associative. 

In particular, $\AA=\CC$ gives the quaternion algebra $\BB=\HH$, 
which is associative, so the Cayley-Dickson process can be applied once more 
and $\AA=\HH$ gives the octonion algebra $\BB=\OO$. This algebra is no longer
associative (although it is alternative, hence normed), and the process
then fails to produce new normed algebras. 

\smallskip
A useful variant of the Cayley-Dickson process is obtained by changing a sign in the 
formula for the product in $\BB$, letting $(x,y)(z,t)=(xz+t\overline{y},\overline{x}t+zy)$. 
The resulting algebra $\tilde\BB$ is called {\em split}.
It is no longer normed,
but endowed with a nondegenerate  quadratic
 form of signature $(a,a)$  compatible 
with the product. (Here
  $a$ denotes the 
dimension of $\AA$.) In particular, it is a {\em composition algebra}. 
The algebra $\tilde\CC$ of split complex numbers is   the algebra 
$\RR\op\RR$ with termwise multiplication, while the split quaternion algebra 
$\tilde\HH$ is equivalent to $M_2(\RR)$.

\subsection{Derivations}
Any automorphism of $\HH$ is inner, so that any derivation of $\HH$ is of the form 
$L_a-R_a$ for some imaginary quaternion $a\in {\rm Im}\HH$, where $L_a$ and $R_a$ denote 
the operators of left and right mulplication by $a$. In particular ${\rm Der}\HH={\rm Im}\HH
=\fso_3$. 
 
\smallskip
The derivation algebra of $\OO$ is the compact Lie algebra $\fg_2$. The algebra ${\rm Der}\HH=\fso_3$ 
does not imbed in ${\rm Der}\OO$ in a canonical way,  but there is one prefered embedding for each 
decomposition $\OO = \HH\oplus e\HH$, 
where the product in $\OO$ is deduced from that in $\HH$ through the Cayley-Dickson
process. Indeed, for $\phi\in Der\HH$, the endomorphism $\tilde\phi$ of $\OO$ defined by 
$\tilde\phi (x+ey) = \phi(x)+e\phi(y)$ is a derivation. Note that the subalgebra 
${\rm Der}(\OO,\HH)$ of 
the derivations of $\OO$ stabilizing $\HH$ is 
strictly greater than ${\rm Der}\HH$. 
It contains, for each 
imaginary quaternion $h$, the map $\psi_h(x+ey)=e(yh)$. This gives another copy of $\fso_3$ in 
${\rm Der}(\OO,\HH)=\fso_3\times\fso_3$, and the restriction of these derivations to $e\HH$ 
gives the full
$\fso_4=\fso_3\times\fso_3$. Explicitly, choose a standard basis $e_0=1$, 
$e_1$, $e_2$, $e_3=e_1e_2$ 
of $\HH$, and let $e_4=e$, $e_5=ee_1$, $e_6=ee_2$, $e_7=ee_3$. 
Using this basis of $\OO$, we obtain a   matrix
representation of $\fg_2={\rm Der}(\OO)$:
$$
\fg_2=\begin{pmatrix}
0 & 0 & 0 & 0 & 0 & 0 & 0 & 0  \\
0 & 0 & -\a_2 &  -\a_3 &  -\a_4 &  -\a_5 &  -\a_6 &  -\a_7 \\
0 & \a_2 & 0 & -\b_3 & -\b_4 & -\b_5 & -\b_6 & -\b_7 \\
0 & \a_3 & \b_3 & 0 & \a_6+\b_5 &  -\a_5-\b_4 &  -\a_4+\b_7 &  \a_7-\b_6 \\
0 & \a_4 & \b_4 &  -\a_6-\b_5 & 0 & -\g_5 & -\g_6 & -\g_7 \\
0 & \a_5 & \b_5 & \a_5+\b_4 & \g_5 & 0 & -\a_2+\g_7  & -\a_3-\g_6 \\
0 & \a_6 & \b_6 & \a_4-\b_7 & \g_6  & \a_2-\g_7 & 0  & \b_3+\g_5 \\
0 & \a_7 & \b_7 & -\a_7+\b_6 & \g_7  & \a_3+\g_6  & -\b_3-\g_5 & 0 
\end{pmatrix} .$$

\medskip
It follows that
the Lie sub-algebra
 ${\rm Der}(\OO,\HH)$   may be written:
$${\rm Der}(\OO,\HH)=\begin{pmatrix}
0 & 0 & 0 & 0 & 0 & 0 & 0 & 0  \\
0 & 0 & -\a_2 &  -\a_3 &  0 &  0 &  0 & 0 \\
0 & \a_2 & 0 & -\b_3 & 0 & 0 & 0 & 0 \\
0 & \a_3 & \b_3 & 0 & 0 &  0 & 0 & 0 \\
0 & 0 & 0 &  0 & 0 & -\g_5 & -\g_6 & -\g_7 \\
0 & 0 & 0 & 0 & \g_5 & 0 & -\a_2+\g_7  & -\a_3-\g_6 \\
0 & 0 & 0 & 0 & \g_6  & \a_2-\g_7 & 0  & \b_3+\g_5 \\
0 & 0 & 0 & 0 & \g_7  & \a_3+\g_6  & -\b_3-\g_5 & 0 
\end{pmatrix} .$$

Our first type of derivations in ${\rm Der}(\OO,\HH)$ corresponds to those matrices for which 
$\g_5=\g_6=\g_7=0$, and the second type to those for which $\a_2=\a_3=\b_3=0$. These 
two copies of $\fso_3$ 
commute. Finally, the south-east corner of the matrix above gives a copy of $\fso_4$. 

\subsection{Triality}

\subsubsection{Cartan's triality}
The Dynkin diagram of type $D_4$ is the only one with a threefold symmetry. 

\setlength{\unitlength}{4mm}
\begin{picture}(30,8)(-12,2)
\put(0,5.7){$D_4$}
\put(6,5.8){$\circ$}
\put(3.6,5.7){$\circ$}
\put(7.8,8.1){$\circ$}
\put(7.8,3.15){$\circ$}
\put(4,6){\line(1,0){2}}
\put(6.4,6.1){\line(2,3){1.5}}
\put(6.4,5.9){\line(2,-3){1.5}}
\end{picture}

Symmetries of Dynkin diagrams detect outer automorphisms of the corresponding 
 complex (or split)
 Lie algebras, and in the case of $D_4$ this is closely related  with the 
algebra structure of the octonions. Indeed, let 
$$T(\OO)=\{(U_1,U_2,U_3)\in SO(\OO)^3,\; U_1(xy)=U_2(x)U_3(y)\;\;\forall x,y\in\OO\}.$$
This group contains ${\rm Aut}(\OO)=G_2$ and has three projections $\pi_1$, $\pi_2$, $\pi_3$ 
on $SO(\OO)$, hence three representations of dimension eight. 
The following theorem is due
 to Elie Cartan and can be found page 370 of his paper \cite{cartan} from 1925.
Define ${\rm Spin}_8$ to be the compact, simply connected Lie group
with Lie algebra $\fso_8$.

\begin{theo} The group $T(\OO)$ is ${\rm Spin}_8$, and each projection 
$\pi_i$ is a twofold covering of $SO_8$. The three corresponding eight dimensional
representations of $\;{\rm Spin}_8$ are not equivalent. \end{theo}

With this description of ${\rm Spin}_8$ in hand, it is easy to define an action of $\cS_3$ by
outer automorphisms. Indeed, since $\OO$ is alternative and $\overline{x}=-x$ modulo the unit 
element of $\OO$, the relation $\overline{x}(xy)=(\overline{x}x)y$ holds. 
(This and similar identities are due to Moufang,
see \cite{harvey}.)
A little computation then 
shows that the identities 
$$U_1(xy)=U_2(x)U_3(y), \quad U_2(xy)=U_1(x)\overline{ U_3(\overline{y})},
 \quad U_3(xy)=\overline{U_2(\overline{x})}U_1(y)$$
are equivalent. If we denote by $\tau\in SO(\OO)$ the conjugation anti-automorphism, 
this implies that we can define automorphisms $s$ and $t$ of ${\rm Spin}_8$ by letting 
$$\begin{array}{rcl}
s(U_1,U_2,U_3) & = & (U_2,U_1,\tau\circ U_3\circ\tau), \\
t(U_1,U_2,U_3) & = & (\tau\circ U_1\circ\tau , \tau\circ U_3\circ\tau ,\tau\circ U_2\circ\tau).
\end{array}$$
These automorphisms are not inner, since they exchange the three inequivalent eight 
dimensional representations of ${\rm Spin}_8$. Moreover, we have $s^2=t^2=1$ and $sts=tst$, 
so that $s$ and $t$ generate a group isomorphic to $\cS_3$, which is the full group
of outer automorphisms of ${\rm Spin}_8$. 

Note that at the level of Lie algebras, Cartan's theorem above implies the 
{\em infinitesimal triality principle}, which we state as
$$\ft(\OO):=\{(u_1,u_2,u_3)\in \fso(\OO)^3,\; u_1(xy)=u_2(x)y+xu_3(y)
\;\forall x,y\in\OO\}=\fso_8.$$
This   means that if we fix $u_1$ for example, there is a unique pair of skew-symmetric 
endomorphism of $\OO$ such that the above relation holds.

\subsubsection{Triality algebras}

Allison and Faulkner \cite{af}, following work
of Allison \cite{al} and Kantor \cite{kantor} define the following vast generalization
of the triality principle.
Let $\cA$ be any algebra with unit and involution $a\mapsto \overline{a}$ 
over a ring of scalars. Let ${\rm Im}\cA=\{a\in\cA,\; \overline{a}=-a\}$.
Define the {\it triality algebra of $\cA$} 
$$
\ft (\cA) := \{ T=(T_1,T_2,T_3)\in \fgl(\cA)^{\op 3}
\mid \overline{T_1}(ab)=T_2(a)b+aT_3(b) \;\;\forall a,b\in\cA\},
$$
with $\overline{T_1}(a)=\overline{T_1(\overline{a})}$. 
This  is a Lie subalgebra of $\fgl (\cA)^{\op 3}$.
In \cite{af} it is call 
the set of {\it partially related Lie triples of $\cA$}.
If $3$ is invertible in the ring of scalars, they provide  (Corollary 3.5) 
a general description of $\ft(\cA)$, which in the case where $\cA$ is 
alternative gives
$$
\begin{aligned}
\ft (\cA)&=   \{(D+L_s-R_t,D+L_t-R_r,D+L_r-R_s) 
     \mid D\in {\rm Der} (\cA ),\ 
r,s,t\in{\rm Im}(\cA ),\ r+s+t=0\}\\
&={\rm Der} (\cA ) \op {\rm Im}(\cA)^{\op 2}. \end{aligned}
$$
  If $D$ is a derivation, $(D,D,D)\in\ft( \cA)$ 
because $\overline{D}=D$. Also, recall that $\cA$ being alternative 
means that any subalgebra generated by two elements is associative, 
or equivalently, that the {\em associator}
$$\{x,y,z\} := (xy)z-x(yz), \quad x,y,z\in\cA,$$
is a skew-symmetric trilinear function. In particular, we have 
$\{x,y,z\}=\{z,x,y\}$, which
can be rewritten as $(xy)z+z(xy)=x(yz)+(zx)y$. If $z$ is imaginary, 
this shows that $(-L_z-R_z,L_z,R_z)$ is in $\ft(\cA)$.
Using the triality automorphisms defined above, one gets that 
$(L_s-R_t,L_t-R_r,L_r-R_s)$ belongs to $\ft(\cA)$ as soon as 
$r,s,t$ are imaginary with $r+s+t=0$. 

\smallskip
For composition algebras, we get 
$$\begin{array}{lcl}
& \ft(\RR)=0 & \\
\ft(\CC)=\RR^2  & & \ft(\SC)=\RR^2 \\
\ft(\HH)=(\fso_3)^{\op 3} &  & \ft(\SH)=(\fsl_2)^{\op 3} \\
\ft(\OO)=\fso_8  & & \ft(\SO)=\fso_{4,4}.
\end{array}$$

\medskip

\subsection{Triality and $\fso_{4,4}$: $4$-ality}
An interesting version of triality holds for $\fso_{4,4}$, 
the split version of $\fso_8$.  We use a model is related to the 
fourfold symmetry of the extended Dynkin diagram ${\tilde D}_4$. 

\setlength{\unitlength}{4mm}
\begin{picture}(30,7.5)(-15,2)
\put(-3,5.7){${\tilde D}_4$}
\put(3.83,5.76){$\circ$}
\put(1.55,8.1){$\circ$}
\put(1.55,3.4){$\circ$}
\put(6.2,8.1){$\circ$}
\put(6.2,3.4){$\circ$}
\put(4.2,6.2){\line(1,1){2}}
\put(4.2,5.8){\line(1,-1){2}}
\put(3.9,6.2){\line(-1,1){2}}
\put(3.9,5.8){\line(-1,-1){2}}
\end{picture}

 Choose four real two-dimensional
vector spaces $A, B, C, D$ and non-degenerate two-forms on each of them.  Then 
$$\fso_{4,4}=\fsl(A)\times \fsl(B)\times \fsl(C)\times \fsl(D)\oplus (A\ot B\ot C\ot D).$$
In this model, $\fsl(A)\times \fsl(B)\times \fsl(C)\times \fsl(D)$ is a sub-Lie algebra 
and acts on $A\ot B\ot C\ot D$ in the obvious way. The bracket of two vectors in 
$A\ot B\ot C\ot D$ can be defined as follows. Choose non-degenerate skew-symmetric
two forms on $A$, $B$, $C$ and $D$. We denote these forms by the same letter $\o$. 
We have then natural isomorphisms $S^2A\ra\fsl(A)$, and so on, by sending a square 
$a^2$, where $a\in A$, to the endomorphism $a'\mapsto \o(a,a')a$. With this understood,
we let 
$$\begin{array}{rcl}
[a\ot b\ot c\ot d, a'\ot b'\ot c'\ot d'] & = & \o(b,b')\o(c,c')\o(d,d')aa' 
+\o(a,a')\o(c,c')\o(d,d')bb' + \\ & & 
+\o(a,a')\o(b,b')\o(d,d')cc' +\o(a,a')\o(b,b')\o(c,c')dd'.
\end{array}$$
One can then check that the Jacobi identity holds, so that we obtain a Lie algebra which
can easily be identified with $\fso_{4,4}$. One could also take in $\fso_{4,4}$, once a 
Cartan subalgebra and a set of positive roots have been chosen, the product of the four 
copies of $\fsl_2$ generated by the three simple roots and by the highest root. These
four copies commute, and the sum of the other root spaces forms a module over their
product. It is then routine to identify this module with the tensor product of the 
natural representations. 

What is particularily nice with this model is that we can explicitly define three 
  inequivalent eight dimensional representations, which we denote by 
$$\begin{array}{rcl}
\SO_1 & = & (A\ot B)\oplus (C\ot D), \\ 
\SO_2 & = & (A\ot C)\oplus (B\ot D), \\
\SO_3 & = & (A\ot D)\oplus (B\ot C).
\end{array}$$
The action of $\fso_{4,4}$ is easy to describe : the action of 
$\fsl(A)\times \fsl(B)\times \fsl(C)\times \fsl(D)$ is the obvious one , 
while that of $A\ot B\ot C\ot D$ on $\SO_1$, for example, is defined by 
$$a\ot b\ot c\ot d.(a'\ot b'+c'\ot d')=\o(c,c')\o(d,d')a\ot b+\o(a,a')\o(b,b')c\ot d.$$
To see that we get three inequivalent representations, we can describe the root 
system of $\fso_{4,4}$ thanks to the previous model, and check that the highest weights 
are distinct. Moreover, once isomorphisms of $A, B, C$ and $D$ with a fixed two-dimensional
vector space have been chosen, there is an obvious action of the symmetric group 
$\cS_4$ on $\fso_{4,4}$ and its triplet of representations. The resulting group 
of outer automorphisms is only $\cS_3$, of course, which in this model appears as the 
quotient of $\cS_4$ by the normal subgroup of order four whose three nontrivial 
elements are products of disjoint two-cycles (recall that the alternating group
$\cA_4$ is not simple!). 
(See \cite{kos} for a geometric interpretation of the $\fso_{4,4}$ triality, 
 and an application to the construction of a minimal 
unitary representation.) Note also that the subalgebra of fixed 
points of the $\cS_3$-action (defined by the permutations of $A,B,C$) on $\fso_{4,4}$ is 
$$\tilde \fg_2  = \fsl(A)\times\fsl(D)\oplus (S^3A\ot D),$$
the split form of $\fg_2$. 

As pointed out by J. Wolf  in \cite{graywolf}, one can twist this $\cS_3$-action so as to get 
a different 
fixed point  subalgebra. We define  automorphisms $\tau$ and $\tau '$ of $\fso_{4,4}$,  
acting on $(M,N,P,Q)\in \fsl(A)\times \fsl(B)\times \fsl(C)\times \fsl(D)$ 
and on $a\ot b\ot c\ot d\in A\ot B\ot C\ot D$ in the following way:
$$\begin{array}{ll}
\tau (M,N,P,Q) = (P,M,N,Q), & \tau(a\ot b\ot c\ot d) = a\ot b\ot c\ot d, \\
\tau '(M,N,P,Q) = (P,M,N,ad(T)Q), & \tau '(a\ot b\ot c\ot d) = c\ot a\ot b\ot T(d), 
\end{array}$$
where $T$ is a rotation matrix  of order three (so that $T^2+T+1=0$). The fixed point set of 
$\tau'$ is $\fsl(A)\times\RR H\oplus M$, where $H=1+2T$ and $M$ is the space of fixed 
points of $\tau'$ inside $A\ot B\ot C\ot D$. Write such a tensor as $X\ot e+Y\ot Te$ for some 
nonzero vector $e$ in $D$, where $X$ and $Y$ belong to $A\ot B\ot C$. It is fixed by 
$\tau'$ if and only if $X=-\tau(Y)$ and $(1+\tau+\tau^2)(Y)=0$, which means that 
$Y=Z-\tau(Z)$ for some $Z$, and then $X=\tau^2(Z)-\tau(Z)$. In particular $M$ has dimension 
$4$ and it is easy to see that the space of fixed points of $\tau'$ is   a copy 
of $\fsl_3$. 

\pagebreak
\subsection{The magic square}

\subsubsection{Jordan algebras}

The relation between exceptional Lie groups and normed division algebras, 
especially the octonion algebra, 
was first noticed  by E. Cartan, at least as early as 1908
(see the  encyclopedia article for the French version of the
Mathematische Wissenschaften).
Cartan observed that
   the 
automorphism group of $\OO$ is the compact group $G_2$. In 1950, Chevalley
and Schafer obtained the compact group $F_4$ as the automorphism group of the 
Jordan algebra 
$$ \cJ_3(\OO) := \Bigg\{\begin{pmatrix} r_1 & x_3 & x_2 \\
 \overline{x_3} & r_2 & x_1 \\
\overline{x_2} & \overline{x_1} & r_3 \end{pmatrix}, 
\quad r_i\in\RR,\; x_j\in\OO\Bigg\}.$$
This is an algebra   with  
 product $A.B=\frac{1}{2}(AB+BA)$. It is a commutative but nonassociative algebra, in which the characteristic identity $A^2(AB)=A(A^2B)$ of 
Jordan algebras   holds. 
By remarks in \cite{car3}, it appears Cartan was aware of this
construction at least as early as 1939.

Jordan algebras were introduced in the 1930's by Jordan, von Neumann and Wigner as 
a natural mathematical framework for quantum theory. Any associative algebra (over a field 
of characteristic not equal to $2$) can be considered as a Jordan algebra for the 
symmetrized matrix product $a.b=\frac{1}{2}(ab+ba)$. A Jordan algebra is said to be {\em exceptional}
if it cannot be embedded as a Jordan subalgebra of an associative algebra. The algebra 
$\cJ_3(\OO)$ is an exceptional Jordan algebra. Moreover,  if $A$ is an 
algebra such that $\cJ_n(A)$, with its symmetrized product, is a Jordan algebra, then $A$ must 
be alternative if $n\ge 3$, and   associative if $n\ge 4$.

The closed orbits $X\subset \cJ_n(\BA)$, with $\BA$ the complexification
of a division algebra are exactly the Scorza varieties
while the Severi varieties are exactly the varieties $X\subset \cJ_3(\BA)$,
see \cite{zak}. 
The close relationship between Scorza varieties and Jordan algebras,
 which is not explained 
in \cite{zak}, is investigated in \cite{chaput2}.

\subsubsection{The Freudenthal-Tits construction}

\smallskip The discovery that $F_4={\rm Aut}\cJ_3(\OO)$ led several authors (including 
Freudenthal, Tits, Vinberg, Jacobson, Springer) to study   the 
relations between exceptional
groups and Jordan algebras. 
Freudenthal and
Tits obtained a unified construction of all compact exceptional
Lie algebras in the following way: take two normed division algebras $\AA$ and $\BB$, 
and let 
$$\fg(\AA,\BB) := {\rm Der}\AA \oplus ({\rm Im}\AA\ot \cJ_3(\BB)_0)\op {\rm Der}\cJ_3(\BB),$$
where $\cJ_3(\BB)_0\subset \cJ_3(\BB)$ is the hyperplane of traceless matrices. 
One can then define a Lie algebra structure on $\fg(\AA,\BB)$, with 
${\rm Der}\AA\op {\rm Der}\cJ_3(\BB)$ as a Lie subalgebra acting on 
${\rm Im}\AA\ot \cJ_3(\BB)_0$ in a natural way \cite{tits0}. 
The result of this construction is 
the {\it Freudenthal-Tits magic square of Lie algebras}:

\medskip
\begin{center}\begin{tabular}{|c|cccc|} \hline 
 & $\RR$ & $\CC$ & $\HH$ & $\OO$ \\ \hline
$\RR$ &  $\fsu_2$ & $\fsu_3$ & $\fsp_6$ & $\ff_4$ \\  
$\CC$ & $\fsu_3$ & $\fsu_3\etimes \fsu_3$ & $\fsu_6$ & $\fe_6$ \\ 
$\HH$ & $\fsp_6$ & $\fsu_6$ & $\fso_{12}$ & $\fe_7$  \\
$\OO$ & $\ff_4$ & $\fe_6$ & $\fe_7$ & $\fe_8$ \\ \hline
\end{tabular}\end{center}
\medskip

Let  $\BB =\underline 0$ denote the zero algebra, so $\cJ_n(\underline 0)\simeq
\BC^n$ is the diagonal matrices. Let
  $\Delta$ be something slightly smaller, so that
 $\cJ_n(\Delta)\subset \cJ_n(\underline 0)$ consists of the scalar
matrices. With these conventions, we obtain a magic rectangle:

\medskip\begin{center}
\begin{tabular}{|c|cccccc|} \hline
 & $\Delta$  &  $\underline 0$  &  $\RR$  &  $\CC$  &  $\HH$  &  $\OO$  \\ \hline
 $\RR$  &  $0$ & $0$ &  $\fsu_2$  &  $\fsu_3$  &  $\fsp_6$  &  $\ff_4 $ \\  
 $\CC$  &  $0$ & $0$ &  $\fsu_3$ &  $\fsu_3\etimes \fsu_3$  &  $\fsu_6$  &  $\fe_6$  \\ 
 $\HH$  &  $\fsu_2$ & $(\fsu_2)^{\op 3}$ &  $\fsp_6$  &  $\fsu_6$  &  $\fso_{12}$  &  $\fe_7$   \\
 $\OO$  &  $\fg_2$ & $\fso_8$ &  $\ff_4$  &  $\fe_6$  &  $\fe_7$  &  $\fe_8$  \\ \hline
\end{tabular}\end{center}
\medskip

\subsubsection{The Allison  construction}

While it is impressive that  the vector space $\fg(\BA,\BB)$ is a Lie
algebra, the symmetry of the square appears to be miraculous.
Vinberg, in 1966,  obtained a symmetric construction
(see \cite{ov} for an exposition).
We present a variant of the Vinberg construction first discovered
by Allison \cite{al}, then rediscovered
independently in more geometric form
by Dadok and Harvey \cite{dh} (who, after
discovering Allison's work decided not  to publish their manuscript),
  Barton and Sudbery \cite{bs}, and ourselves  \cite{LMtrial}.

Allison shows that given an arbitrary structurable algebra $\cA$,
one can put a Lie algebra structure on 
$$\ft (\cA )\op \cA^{\op 3}$$  
  Allison 
observes moreover that the algebra $\cA=\BA\ot \BB$ with $\BA,\BB$ division
algebras is structurable and
his construction applied to them
 yields the triality model, see \cite{al, af}.
 \smallskip
 
Using the Allison construction, the magic square may
be described as follows:
$$\fg(\AA,\BB) = \ft(\AA)\times\ft(\BB)\op(\AA_1\ot\BB_1)\op(\AA_2\ot\BB_2)
 \op(\AA_3\ot\BB_3).$$

This bracket is defined so that $\ft(\AA)\times\ft(\BB)$ is a Lie subalgebra, 
acting on each $\AA_i\ot\BB_i$ in the natural way. The bracket of an element 
$a_1\ot b_1\in\AA_1\ot\BB_1$ with $a_2\ot b_2\in\AA_2\ot\BB_2$ is simply
$a_1a_2\ot b_1b_2$, considered as an element of $\AA_3\ot\BB_3$. This is the 
general rule for taking the bracket of an element of $\AA_i\ot\BB_i$ with one 
of $\AA_j\ot\BB_j$, although there are
 some slight twists whose details can be found in
 \cite{af} or \cite{LMtrial}. 
Finally, the bracket of two elements in $\AA_i\ot\BB_i$ is defined by the quadratic 
forms on $\AA_i$, $\BB_i$ and the natural maps $\Psi_i: \wedge^2\AA_i\ra\ft(\AA)$ obtained 
by dualizing the action of $\ft(\AA)$ on $\AA_i$, which can be seen as  a map 
$\ft(\AA)\ra\wedge^2\AA_i$. (Note that $\ft(\AA)$ is always reductive, hence 
isomorphic to its dual as a $\ft(\AA)$-module.) The key formulas that ensure that 
the Jacobi identities hold in $\fg(\AA,\BB)$ are then the following, which we state 
only for $\Psi_1$:
$$\begin{array}{rcl}
\Psi_1(u\wedge v)_1x & = & Q(u,x)v-Q(v,x)u, \\
\Psi_1(u\wedge v)_2x & = & \frac{1}{2}(\overline v (ux)-\overline u (vx)), \\
\Psi_1(u\wedge v)_3x & = & \frac{1}{2}((xu)\overline v -(xv)\overline u ).
\end{array}$$
These formulas are classical in the study of triality, that is for $\AA=\OO$. 
It is easy to check that in the other cases, which are much easier,  we can
arrange so that they also hold. 

\begin{prop} 
The real Lie algebra $\fg(\AA,\BB)$ is compact, while $\fg(\SA,\SB)$ is split.
\end{prop}

\proof Let $K_{\AA}$ be the non-degenerate quadratic form on $\ft(\AA)$ which we used in the 
construction of $\fg(\AA,\BB)$. For $\AA=\OO$, this form is the Killing form on $\fso_8$ and is 
therefore negative definite; this is also true for the other division algebras. 
Now it is easy to 
check that the quadratic form 
$$\cQ := K_{\AA}+K_{\BB}-\sum_{i=1}^3Q_{\AA_i}\ot Q_{\BB_i}$$
is invariant and negative definite. Since $\fg(\AA,\BB)$ is simple (except for $\AA=\BB=\CC$),
this invariant form must be proportional to the Killing form, which is therefore 
negative definite 
as well (it cannot be positive definite!). Thus  $\fg(\AA,\BB)$ is compact. 

Finally $\fg(\SA,\SB)$ is clearly split since both $\ft(\SA)$ and $\ft(\SB)$ are split,
and their product is a maximal rank subalgebra of $\fg(\AA,\BB)$. \qed 

\subsection{Inclusions}
It is not at all 
obvious that an inclusion $\BB\subset\BB '$ of composition 
algebras induces an inclusion $\fg(\AA,\BB)\subset\fg(\AA,\BB ')$ 
of Lie algebras. In fact the Allison-Faulkner process
applied to structurable algebras $\cA \subset \cB$ does
{\it not} imply inclusions $\fg (\cA)\subset \fg (\cB)$. 
To discuss this problem in our
situation we may assume that $\BB '=\BB\op e\BB$ 
is deduced from $\BB$ 
via the Cayley-Dickson process. 

\begin{theorem} There is a unique
 embedding of Lie algebras $\ft(\BB)\hookrightarrow\ft(\BB ')$ 
that makes $\fg(\AA,\BB)$ a Lie subalgebra of $\fg(\AA,\BB ')$.

If $\BB '$ is deduced from $\BB$ by the Cayley-Dickson process,
 the relevant 
embedding of $\ft(\BB)$ inside $\ft(\BB ')$ is  
$$\ft(\BB)\simeq\{(U_1,U_2,U_3)\in\ft(\BB ')\mid 
\quad U_i(\BB)\subset\BB\}\subset\ft(\BB').$$
 \end{theorem}

\proof 
We want the inclusion $\fg(\AA,\BB)\hookrightarrow\fg(\AA,\BB ')$ to be a
morphism of Lie algebras. 
It is easy to see that the only  non-obvious verification to be made is that 
the bracket of two elements in $\AA_i\ot\BB_i$ is the same in $\fg(\AA,\BB)$ as in  
$\fg(\AA,\BB ')$, which amounts to proving that we have a commutative diagram 
$$\begin{array}{ccc}
\wedge^2(\AA_i\ot\BB_i) & \lra & \ft(\AA)\times\ft(\BB) \\
\downarrow & & \downarrow \\
\wedge^2(\AA_i\ot\BB_i') & \lra & \ft(\AA)\times\ft(\BB ') . \\
\end{array}$$
Since the norm on $\BB$ is the restriction of the norm on $\BB '$, it is actually enough 
to prove the commutativity of the following diagram:
$$\begin{array}{ccc}
\wedge^2\BB_i & \stackrel{\Psi_i}{\lra} & \ft(\BB) \\
\downarrow & & \downarrow \\
\wedge^2\BB_i' & \stackrel{\Psi_i'}{\lra} & \ft(\BB ') .
\end{array}$$

We must check this for all $i$, so we consider the sum $\Psi_{\BB} : \wedge^2\BB_1\op
\wedge^2\BB_2\op \wedge^2\BB_3\lra\ft(\BB)$. Note that $\Psi_{\BB}$ is dual to the action, 
hence surjective since the action map is injective by definition. So the commutativity 
of the diagram above for $i=1,2,3$  completely determines
the inclusion of $\ft(\BB)$ 
inside $\ft(\BB ')$. For this inclusion map to be well-defined, we need a compatibility
condition, namely that ${\rm Ker}\Psi_{\BB} \subset {\rm Ker}\Psi_{\BB '}$. 
Everything being trivial for $\BB=\RR$, we will check that 
${\rm Ker}\Psi_{\CC} \subset {\rm Ker}\Psi_{\HH}\subset {\rm Ker}\Psi_{\OO}$
and deduce the relevant embeddings $\ft(\CC)\subset\ft(\HH)\subset\ft(\OO)$. 

\smallskip
We begin with $\Psi_{\CC}$, which is dual to the inclusion map $$\ft(\CC)=\{(a,b,c)\in\RR^3,
\;a=b+c\}\subset \wedge^2\CC_1\op\wedge^2\CC_2\op \wedge^2\CC_3,$$ where for example $a\in\RR$ 
maps to the skew-symmetric endomorphism {\tiny $\begin{pmatrix} 0 & a \\ -a & 0 
\end{pmatrix}$} of $\CC_1$. 
On $\ft(\CC)$ we use the restriction of the canonical norm of $\RR^3$, and we find that the dual 
map is 
$$\Psi_{\CC} : \Bigg (\begin{pmatrix} 0 & a \\ -a & 0 \end{pmatrix}, 
\begin{pmatrix} 0 & b \\ -b & 0 \end{pmatrix}, \begin{pmatrix} 0 & c \\ -c & 0 \end{pmatrix}
\Bigg) \mapsto \Big(\frac{2a+b+c}{3},\frac{a+2b-c}{3},\frac{a-b+2c}{3}\Big),$$
whose kernel is generated by {\tiny $\Bigg(\begin{pmatrix} 0 & 1 \\ -1 & 0 \end{pmatrix}, 
\begin{pmatrix} 0 & -1 \\ 1 & 0 \end{pmatrix}, \begin{pmatrix} 0 & -1 \\ 1 & 0 
\end{pmatrix}\Bigg)$} .

\medskip
Now to compute $\Psi_{\HH}$,  using that
${\rm Der}(\HH)\simeq {\rm Im}\HH$ we may describe $\ft(\HH)$ as 
$$\ft(\HH) = \{(L_a+R_b, L_a-R_c, L_c+R_b), \quad a,b,c\in {\rm Im}\HH\}.$$
Every skew-symmetric endomorphism of $\HH$ is of the form $L_x+R_y$ for some 
imaginary quaternions $x$ and $y$, and 
the norm defined by the Killing form is   
$\norm{L_x+R_y}^2 = \norm{x}^2 +\norm{y}^2$. With the induced norm on $\ft(\HH)$ 
we   compute that 
$$\Psi_{\HH} : (L_{x_1}+R_{y_1}, L_{x_2}+R_{y_2}, L_{x_3}+R_{y_3})
\mapsto (L_{x_1+x_3}+R_{y_1+y_2}, L_{x_1+x_3}-R_{x_2-y_3}, L_{x_2-y_3}+R_{y_1+y_2}).$$
Its kernel is $\{(L_x+R_y, L_z-R_y, -L_x+R_z), \quad x,y,z\in {\rm Im}\HH\}.$

\medskip
Finally, the computation of $\Psi_{\OO}$ is an immediate consequence of the triality 
principle. Indeed, we know that each projection $\ft(\OO)\ra\wedge^2\OO_i$ is an isomorphism. 
If we denote by $(U,U',U'')$ the elements of $\ft(\OO)$, this means that such a triplet is 
uniquely determined by $U$, $U'$ or $U''$. We deduce that  
$$\Psi_{\OO} : (U_1, U_2', U_3'')\mapsto (U_1+U_2+U_3, U_1'+U_2'+U_3', U_1''+U_2''+U_3'').$$

We are ready to check our compatibility conditions. First note that the endomorphism
{\tiny $\begin{pmatrix} 0 & 1 \\ -1 & 0 \end{pmatrix}$} of $\CC$, when extended to $\HH$ in the 
obvious way, is   $1\wedge i=\frac{1}{2}(L_i+R_i)$. Taking $x=y=-z=i$ in our description
of ${\rm Ker}\Psi_{\HH}$, we see that it contains ${\rm Ker}\Psi_{\CC}$, as required, and that 
the induced inclusion $\ft(\CC)\subset\ft(\HH)$ is given by
$$(b,c)\in\RR^2\mapsto (L_{(b+2c)i}+R_{(2b+c)i},  L_{(b+2c)i}+R_{(c-b)i},  
L_{(b-c)i}+R_{(2b+c)i}).$$

Now take  $(L_x+R_y, L_z-R_y, -L_x+R_z)$ in ${\rm Ker}\Psi_{\HH}$, where 
$x,y,z\in {\rm Im}\HH$, and consider it as an element of $\wedge^2\OO_1\op
\wedge^2\OO_2\op\wedge^2\OO_3$, which we denote by $((L_x+R_y,0), 
(L_z-R_y,0), (-L_x+R_z,0))$. 
Here   $(L_z,0)$, for example,
denotes  the endomorphism of $\OO$ defined on elements of 
$\HH$ by left multiplication by $z\in\HH$, and by zero on $\HH^{\perp}$. 
Note that with these 
conventions, it follows from the multiplication rules in the Cayley-Dickson process that for 
$z\in {\rm Im}\HH$, we have $L_z=(L_z,-L_z)$ and $R_z=(R_z,L_z)$. 

It is easy to check that every skew-symmetric endomorphism of $\HH$ is of the form $L_a+R_b$. 
In matrix terms, if $a=\a_1i+\a_2j+\a_3k$ and $b=\b_1i+\b_2j+\b_3k$, we have 
$$L_a=\begin{pmatrix} 0 & -\a_1& -\a_2& -\a_3 \\ \a_1 & 0 & -\a_3  & \a_2 \\
\a_2 & \a_3 & 0 & -\a_1 \\ \a_3 & -\a_2 & \a_1 & 0 \end{pmatrix}, \quad
R_b=\begin{pmatrix} 0 & -\b_1& -\b_2& -\b_3 \\ \b_1 & 0 & \b_3  & \b_2 \\
\b_2 & -\b_3 & 0 & \b_1 \\ \b_3 & \b_2 & -\b_1 & 0 \end{pmatrix}.$$
Under the natural isomorphism $\fso(\HH)\simeq\wedge^2\HH$ induced by the norm on $\HH$, 
this means that  
$$\begin{array}{rcl}
L_a & = & \a_1(1\we i+j\we k)+  \a_2(1\we j+k\we i)+  \a_3(1\we k+i\we j), \\  
R_b & = & \b_1(1\we i-j\we k)+  \b_2(1\we j-k\we i)+  \b_3(1\we k-i\we j).
\end{array}$$
Now we can ask what is the element of $\ft(\OO)$ defined by 
$u_1=1\we i\in\fso(\HH)\subset\fso(\OO)$? To answer this question, 
we just need to note that $1\we i=\frac{1}{2}(L_i+R_i)$, so 
$u_2=\frac{1}{2}L_i$ and $u_3=\frac{1}{2}R_i$. Also if we 
take $u_1=j\we k$, a straightforward computation gives  
$u_2=j\we k+\frac{1}{2}R_i=\frac{1}{2}(L_i,L_i)$ and 
$u_3=j\we k-\frac{1}{2}L_i=\frac{1}{2}(-R_i,L_i)$. 
Letting automorphisms of $\HH$ act, we deduce that 
the following triplets $(u_1,u_2,u_3)$ 
belong to $\ft(\OO)$, where $z,u,v$ denote imaginary quaternions:
$$\begin{array}{cccccc}
u_1\hspace*{1cm} & 1\we z & u\we v & (L_z,0) & (R_z,0) & (0,L_z)\\
u_2\hspace*{1cm} & \frac{1}{2}L_z & u\we v+\frac{1}{4}R_{[u,v]} & 
(L_z,0) & (0,-L_z) &(R_z,0) \\
u_3\hspace*{1cm} & \frac{1}{2}R_z & u\we v-\frac{1}{4}L_{[u,v]} & 
(0,L_z) & (R_z,0) &(-L_z,0) 
\end{array}$$

Thus we have 

$$\begin{array}{ccccc}
\Psi_{\OO}((L_x,0),0,0) & = &  ((L_x,0), 
(0,-L_x), (L_x,0)) & = & \Psi_{\OO}(0,0,(L_x,0)), \\
\Psi_{\OO}((R_y,0),0,0) & = &  ((R_y,0), 
(R_y,0), (0,L_y)) & = & \Psi_{\OO}(0,(R_y,0),0), \\
\Psi_{\OO}(0,(L_z,0),0) & = &  ((0,-L_z), (L_z,0), 
(-R_z,0)) & = & \Psi_{\OO}(0,0,(-R_z,0)).
\end{array}$$

\smallskip
This immediately implies that $((L_x+R_y,0), 
(L_z-R_y,0), (-L_x+R_z,0))$ belongs to 
${\rm Ker}\Psi_{\OO}$, and we deduce from the preceding formulas that the induced 
imbedding of $\ft(\HH)\subset\ft(\OO)$ is given by 
$$\begin{array}{l}
(L_a+R_b, L_a-R_c, L_c+R_b)\in\ft(\HH) \\ \hspace*{3cm}\mapsto 
((L_a+R_b,L_c), (L_a-R_c,-L_b), (L_c+R_b,L_a))\in\ft(\OO).
\end{array}$$
In particular, one can check that as a subalgebra of $\ft(\OO)$,  $\ft(\CC)$ is generated 
by $(L_i+R_i,L_i,R_i)$ and $(R_i,-R_i,L_i+R_i)$. 

\smallskip
A   more careful inspection of these formulas yields
the explicit conclusion stated in the theorem.\qed
 
\subsection{Automorphisms, symmetric and trisymmetric spaces}

The 
 symmetry in the triality model in the roles
 of $\BA$ and $\BB$  allows one   to exhibit automorphisms of exceptional
Lie algebras with interesting properties. For example:

\begin{prop} \cite{LMtrial} The endomorphism of $\fg(\AA,\BB)$ defined as the identity on 
$\fh(\AA,\BB)=\ft(\AA)\times\ft(\BB)\oplus \AA_1\ot\BB_1$, and minus the identity 
on $\AA_2\ot\BB_2\oplus \AA_3\ot\BB_3$, is a Lie algebra involution.
\end{prop}

Exhibiting a Lie algebra involution or a symmetric space is more or less equivalent,
and we conclude that there exists a ``magic square'' of symmetric spaces of 
dimensions $2ab$  where $a$ and $b$ denote the dimensions of $\AA$ and $\BB$.
In 
particular we always get powers of $2$! These are not all 
uniquely defined by our Lie algebra involution, but a version of this magic square 
 is the following: 

\medskip
\begin{center}\begin{tabular}{|c|cccc|} \hline 
 & $\RR$ & $\CC$ & $\HH$ & $\OO$ \\ \hline
$\RR$ & $\RR\PP^2$ & $\CC\PP^2$ & $\HH\PP^2$ & $\OO\PP^2$ \\  
$\CC$ & $\CC\PP^2$ & $\CC\PP^2\times\CC\PP^2$ & $G_{\CC}(2,6)$ & $\OO\PP^2_{\CC}$ \\ 
$\HH$ & $\HH\PP^2$ & $G_{\CC}(2,6)$ & $G_{\RR}(4,12)$ & $E_{7(-5)}$ \\
$\OO$ & $\OO\PP^2$ & $\OO\PP^2_{\CC}$ &  $E_{7(-5)}$ 
& $E_{8(8)}$ \\ \hline
\end{tabular}\end{center}
\medskip
 The symmetric space $E_{8(8)}=
E_8/SO_{16}$, of dimension $128$, is particularly intriguing and it would be very nice 
to have a direct construction of it. It is claimed in \cite{rosenfeld} that it can be 
interpreted as a projective plane over $\OO\ot\OO$, and that in fact the whole square
above of symmetric spaces can be obtained by taking projective planes, in a suitable sense, 
over the tensor products $\AA\ot\BB$. 

\medskip
It is also easy to describe Lie algebra automorphisms of order three, reflecting 
the trialitarian origin of the magic square. Recall that $\ft(\AA)$ has a natural 
automorphism $\tau_{\AA}$ or order three. 

\begin{prop} The endomorphism of $\fg(\AA,\BB)$ defined by 
$$\tau_{\AA,\BB}(X+Y+U_1+U_2+U_3)=\tau_{\AA}(X)+\tau_{\BB}(Y)+U_2+U_3+U_1$$
for $X\in\ft(\AA)$, $Y\in\ft(\BB)$, $U_i\in\AA_i\ot\BB_i$, is a Lie algebra automorphism
of order three, whose fixed points set is the subalgebra 
$$\fk(\AA,\BB)=Der\AA\times Der\BB\op \AA\ot\BB.$$ 
When $\BB=\OO$, the   conclusion
also holds when we replace the usual triality 
automorphism $\tau_{\OO}$ by its twisted version $\tau'_{\OO}$.
In this case,   the fixed point set is the subalgebra 
$\fk'(\AA,\BB)=Der\AA\times\fsu_3\op \AA\ot\BB$. 
\end{prop}

As in the case of involutions, we deduce from this statement a ``magic square'' of 
homogeneous spaces which are no longer symmetric, but are sometimes called trisymmetric:
 they are quotients of a semisimple Lie group $G$ by a subgroup $K$ which is, up to a finite 
group, the fixed point set of an automorphism of $G$ of order three. 
Their dimensions are equal to $2ab+2a+2b-4$.

{\small
\smallskip\begin{center}\begin{tabular}{|c|cccc|}\hline
 & $\RR$ & $\CC$ & $\HH$ & $\OO$ \\ \hline
$\RR$ & $U_2/U_1\times U_1$ & $U_3/(U_1)^3$ & $Sp_6/U_2\times Sp_2$ & 
$F_4/(Spin_7\times S^1/\ZZ_2)$ \\  
$\CC$ & $U_3/(U_1)^3$ & $U_3\times U_3/(U_1)^6$ & $U_6/(U_2)^3$ & 
$E_6/SO_8\times (SO_2)^2$ \\ 
$\HH$ & $Sp_6/U_2\times Sp_2$ & $U_6/(U_2)^3$ & 
$SO_{12}/U_4\times SO_4$ & $E_7/S(U_7\times U_1)/\ZZ_2$ \\
$\OO$ & $F_4/(Spin_7\times S^1/\ZZ_2)$ & $E_6/SO_8\times (SO_2)^2$ & $E_7/S(U_7\times U_1)/\ZZ_2$
& $E_8/SO_{14}\times SO_2$ \\ \hline 
\end{tabular}\end{center}}
\smallskip

Here we denoted by $F_4, E_6, E_7, E_8$ the compact centerless
groups of these types. Note that for example, $U_3/U_1\times U_1\times U_1$ 
is   the 
variety of complete flags in $\CC^3$. Similarly, $U_6/U_2\times U_2\times U_2$ is the flag
variety $\FF(2,4,6)$. Also $SO_{12}/U_4\times SO_4$ is the space of $8$-dimensional subspaces
in $\RR^{12}$ endowed with an orthogonal complex structure.  
This is because $SO_{2n}$ acts transitively on the space of orthogonal 
complex structures. Moreover, the stabilizer of a point is 
the group of orthogonal transformations commuting with the corresponding
complex structure, and this subgroup of $SO_{2n}$ is a copy of $U_n$.

\smallskip
In the twisted case, for $\BB=\OO$, we replace the last line by 
\smallskip
\begin{center}\begin{tabular}{cccclcccc}
$F_4/(SU_3\times SU_3/\ZZ_3)$ & $E_6/(SU_3\times SU_3\times SU_3/\ZZ_3)$ &  
$E_7/(SU_3\times SU_6/\ZZ_3)$ & $E_8/SU_3\times E_6$, 
\end{tabular}\end{center}
\smallskip
\noindent a series of trisymmetric spaces of dimension $18a+18$. Note 
that at the infinitesimal level, we obtain the quotients $\fg(\AA,\OO)/\fsu_3\times\fg(\AA,\CC)$. 
 
Automorphisms of order three were studied and classified by   Wolf and Gray \cite{graywolf}. 
In the exceptional cases the construction above is not far from  giving the whole classification. 

\subsection{Dual pairs}
The triality model allows one 
to identify series of {\em dual pairs} in the exceptional Lie algebras. 
Recall that a pair $(\fh,\fh')$ of Lie subalgebras of a Lie algebra $\fg$ is a dual pair 
if $\fh$ is the centralizer of $\fh '$ and vice versa. Dual pairs of Lie algebras, or of 
Lie groups (with the same definition), have been  
extensively studied after the discovery by R. Howe 
that certain dual pairs of groups inside the real symplectic group  $Sp(2n,\RR)$ have very 
special properties with respect to the (infinite dimensional) {\em metaplectic representation}
\cite{howe}. 

First recall that, 
as we saw in 4.7, the triality model is compatible with inclusions, in the sense that 
$\fg(\AA,\BB)$ is naturally embedded in $\fg(\AA',\BB')$ when $\AA\subset\AA'$ and 
$\BB\subset\BB'$. 

\begin{prop} The centralizers in $\fg(\AA,\OO)$ of the subalgebras $\fg(\AA,\HH)$, 
$\fg(\AA,\CC)$, $\fg(\AA,\RR)$ and $\ft(\AA)$,  are isomorphic to $\fsu_2$, 
$\fsu_3$, $\fg_2$ and $\fz(\ft(\AA))\times\fso_8$, respectively. Moreover, the centralizers 
of these centralizers are the subalgebras themselves.  
\end{prop}

Here we denoted by $\fz(\ft(\AA))$ the center of $\ft(\AA)$, which is $\ft(\AA)$ itself
when this algebra is commutative, that is for $\AA=\RR$ or $\CC$, and zero otherwise. 
\medskip

\proof It is easy to see that the centralizer of $\fg(\AA,\BB)$ inside $\fg(\AA,\OO)$ 
is a subalgebra of $\ft(\OO)$, equal to the subalgebra of the centralizer of $\ft(\BB)$ 
inside $\ft(\OO)$ acting trivially on each $\BB_i\subset\OO_i$. 

\smallskip
If $\BB=\RR$ we get the subalgebra of elements $(u_1,u_2,u_3)\in\ft(\OO)$ killing the unit 
element in each $\OO_i$. Making $x=1$ and $y=1$ in the triality relation $u_1(xy)=u_2(x)y
+xu_3(y)$, we see that $u_1=u_2=u_3$ is a derivation of $\OO$, so that the center of 
$\fg(\AA,\RR)$ inside $\fg(\AA,\OO)$ is ${\rm Der}\OO=\fg_2$. 

\smallskip
When $\BB=\CC$, the condition that $(u_1,u_2,u_2)$ kills $\CC_i\subset\OO_i$ implies
that $u=u_1=u_2=u_3$ is a derivation of $\OO$ such that $u(i)=0$ ($i=e_1$ in the notations 
of 5.1). But then $u(ix)=
iu(x)$, which means that $u$ is a complex endomorphism of $\CC^{\perp}\simeq\CC^3$, 
endowed with the complex structure defined by $L_i$ (and also that defined by $R_i$, which is
the conjugate complex structure, since $u(xi)=u(x)i$). Note that the corresponding 
element of $\ft(\OO)$ automatically centralizes $\ft(\CC)$, since $u$ commutes with 
$L_i,R_i$ and we saw that $\ft(\CC)$ is generated as a subalgebra of $\ft(\OO)$ by 
$(L_i+R_i,L_i,R_i)$ and $(R_i,-R_i,L_i+R_i)$. Finally, take the general matrix of 
$\fg_2$ as  written in \S 3.2. The condition that $u(i)=0$ means that $\a_2=
\a_3=\a_4=\a_5=\a_6=\a_7=0$. A complex basis of the space $\CC^{\perp}$ is 
$(e_2,e_4,e_6)$, in terms of which $u$ is the complex endomorphism given by the matrix 
$$\begin{pmatrix} i\b_3 & -\b_4+i\b_5 & -\b_6+i\b_7\\
\b_4+i\b_5 & i\g_5 & -\g_6+i\g_7\\ \b_6+i\b_7 & \g_6+i\g_7 & -i(\b_3+\g_5)\end{pmatrix},$$
which is a general element of $\fsu_3$. 

\smallskip
Finally let $\BB=\HH$. Again, our description of $\fg_2$ in \S 3.2 shows that the 
derivations of $\OO$ vanishing on $\HH$ form a copy of $\fso_3$, 
which we can describe as the space of endomorphisms of the form $(0,R_z)$ with 
$z\in {\rm Im}\HH$. But we saw that $\ft(\HH)$, as a subalgebra of $\ft(\CC)$, 
was generated by the triplets 
$((L_a+R_b,L_c), (L_a-R_c,-L_b), (L_c+R_b,L_a))$, where $a,b,c\in {\rm Im}\HH$. 
But these elements of $\ft(\OO)$   commute with those of the form 
$((0,R_z),(0,R_z),(0,R_z))$. This proves that the centralizer of $\fg(\AA,\HH)$ 
in $\fg(\AA,\OO)$ is $\fso_3$. 

Alternatively we can identify this centralizer to the subalgebra of the centralizer 
of $\ft(\CC)$ annihilating $e_2=j$, which is clearly $\fsu_2$. 

The rest of the claim is straightforward. 
\qed

\medskip
In particular, we have pairs of reductive subalgebras in $\fg(\AA,\OO)$, each of which
is the centralizer of the other: this is a {\sl reductive dual pair}. 
The classification of reductive dual 
pairs inside reductive complex Lie algebras was obtained in \cite{ruben}. We deduce from
the proposition above the existence of uniform series of what Rubenthaler calls {\sl towers} of 
reductive dual pairs. 

$$\begin{array}{ccc}
\fg(\AA,\OO) & \hspace*{5mm} = \hspace*{5mm}  & \fg(\AA,\OO) \\ 
\cup &  &\cup \\
\fg(\AA,\HH) & & \fso_8\times\fz(\ft(\AA)) \\ \cup &  &\cup \\
\fg(\AA,\CC) & & \fg_2\\ \cup &  &\cup \\
\fg(\AA,\RR) &  & \fsu_3\\ \cup &  &\cup \\
\ft(\AA) &  & \fsu_2
\end{array}$$

\medskip
We conclude that $\fso_8\times\ft(\AA)$, $\fg_2\times\fg(\AA,\RR)$, 
$\fsu_3\times\fg(\AA,\CC)$ and $\fso_3\times\fg(\AA,\HH)$ are maximal
rank reductive subalgebras of $\fg(\AA,\OO)$. A natural question to ask is 
to describe the module structure of $\fg(\AA,\OO)$ over each of these 
subalgebras. For the first one, the answer is given by the triality
 construction. For the second one, we note that 
$\fg(\AA,\RR)={\rm Der}(\cJ_3(\AA))$ and $\fg_2 = {\rm Der}\OO$, so that we are back to
the original Tits construction of the exceptional series of Lie algebras. 
We obtain:
$$\begin{array}{rcrcl}
\fg(\AA,\OO) & = & \fg_2\times \fg(\AA,\RR) & \op & \RR^7\ot \cJ_3(\AA)_0,  \\
 & = & \fsu_3\times\fg(\AA,\CC) & \op &  \RR^3\ot \cJ_3(\AA)  \op  
(\RR^3\ot \cJ_3(\AA))^*, \\
 & = & \fsu_2\times\fg(\AA,\HH) & \op &  \RR^2\ot \cZ_3(\AA).
\end{array}$$
Here $\cZ_3(\AA)$ denotes the space of  {\em Zorn matrices}, which is a $\fg(\AA,\HH)$-module,
see  \S 4.6. Note that $\fsu_2$ acts not on $\RR^2$ but on $\CC^2$. Therefore 
 there must exist a complex structure on $\cZ_3(\AA)=\cZ_3^0(\AA)\op i\cZ_3^0(\AA)$, 
which is not $\fg(\AA,\HH)$-invariant, so that $\RR^2\ot \cZ_3(\AA)=\CC^2\ot \cZ_3^0(\AA)$,
the $\fsu_2$-action coming from $\CC^2$.

 These towers of dual pairs in series were discovered in
 joint work with B. Westbury.

\subsection{The quaternionic form}
The algebras $\fg(\AA,\SB)$ are also interesting, especially when $\BB=\OO$. 
Following Wolf \cite{wo1}, each simple complex Lie algebra has a unique real form 
such that the associated compact or non-compact symmetric spaces have an invariant quaternionic
structure. We will say that this real form is {\it quaternionic}. 

Recall that the {\it rank} of a symmetric space $G/K$, corresponding to some Cartan decomposition 
$\fg =\fk\op\fp$, is the maximal dimension of a subalgebra of $\fg$ contained in $\fp$. 
Such a subalgebra is automatically abelian and is called a Cartan subspace. 

\begin{theo}
The real Lie algebra $\fg(\AA,\SO)$ is quaternionic. A Cartan subspace of $\ft(\SO)=\fso_{4,4}$
embeds as a Cartan subspace of $\fg(\AA,\SO)$, which in particular has real rank four
(independent  of $\AA$). 
\end{theo}

\proof Recall that each simple complex Lie algebra $\fg$, once a
 Cartan subalgebra and a set of positive 
roots has been chosen, has a canonical 5-graduation defined by the highest 
root $\ta$, more precisely
by the eigenspaces of $ad(H_{\ta})$, where $H_{\ta}$ denotes the coroot of $\a$. 

If $\fg$ is the complexification of the compact real Lie algebra $\fg(\AA,\OO)$, we can 
arrange so that the highest root comes from $\ft(\OO)_{\CC}=\fso_8(\CC)$. 
Note that the space of highest root vectors in $\fso_8(\CC)$ 
is the cone over the  Grassmannian of isotropic planes in $\CC^8$. Using this, we then 
 check that we can 
arrange so that 
$$iH_{\ta} = \begin{pmatrix} 
0 & 0 & 0  & 0  & 0  & 0  & 0  & 0 \\
0  & 0  & 0  & 0  & 0  & 0  & 0  & 0 \\
0  & 0  & 0 & 0  & 0  & 0  & 0  & 0 \\
0  & 0  & 0 & 0  & 0  & 0  & 0  & 0  \\
0 & 0  & 0  & 0  & 0  & -1  & 0  & 0  \\
0 & 0  & 0  & 0  & 1  & 0  & 0  & 0  \\
0 & 0  & 0  & 0  & 0  & 0  & 0  & -1  \\
0 & 0  & 0  & 0  & 0  & 0  & 1  & 0  \end{pmatrix} .$$

The eigenvalues of $ad(H_{\ta})$ acting on $\fg$ are $0,\pm 1, \pm 2$. 
The corresponding eigenspaces are 
$$\begin{array}{rcl}
\fg_0 & = & \ft(\AA_{\CC})\times\ft(\OO_{\CC})_0 \oplus \AA_{\CC,1}\otimes\OO_{\CC,1}^0
\oplus \AA_{\CC,2}\otimes\OO_{\CC,2}^0\oplus \AA_{\CC,3}\otimes\OO_{\CC,3}^0, \\
\fg_{\pm 1} & = &  \ft(\OO_{\CC})_{\pm 1} \oplus \AA_{\CC,1}\otimes\OO_{\CC,1}^{\pm 1}
\oplus \AA_{\CC,2}\otimes\OO_{\CC,2}^{\pm 1}\oplus \AA_{\CC,3}\otimes\OO_{\CC,3}^{\pm 1}, \\
\fg_{\pm 2} & = &  \ft(\OO_{\CC})_{\pm 2},
\end{array}$$
where we denoted by $\OO_i^t$ the $t$-eigenspace for the action of $H_{\ta}$ on $\OO_i$. 
Note that $iH_{\ta}$ belongs to ${\rm Der}\OO =\fg_2$ (see our description of $\fg_2$ above),
so we can let $i=1$ and consider the standard action of $H_{\ta}$. Then it is clear that 
the kernel of $H_{\ta}$ is   the complexification of the standard quaternion subalgebra 
$\HH$, and that the sum of the two other eigenspaces is the complexification of $e\HH$. 

Now, let $\theta$ denote the involution of $\fg$ associated to
 its $5$-graduation. By 
definition, $\theta$ acts by $(-1)^k$ on $\fg_k$. Then $\theta$ stabilizes $\fg(\AA,\OO)$, 
and the corresponding Cartan decomposition $\fg(\AA,\OO)=\fk\oplus\fp$ into eigenspaces 
of $\theta$ is given, in terms of the Cartan decomposition 
$\ft(\OO)=\fk_{\OO}\oplus\fp_{\OO}$, by
$$\begin{array}{rcl}
\fk & = & \ft(\AA)\times\fk_{\OO}\oplus \AA_1\otimes\HH_1\oplus \AA_2\otimes\HH_2
\oplus \AA_3\otimes\HH_3, \\
\fp & = & \fp_{\OO}\oplus \AA_1\otimes e\HH_1\oplus \AA_2\otimes e\HH_2
\oplus \AA_3\otimes e\HH_3.
\end{array}$$

\smallskip
Finally (see \cite{wo1,clerc}), we obtain the quaternionic form 
$\fg(\AA,\OO)_{\HH}$ of $\fg$ by 
twisting this decomposition, 
that is, letting $\fg(\AA,\OO)_{\HH}=\fk\oplus i\fp$, which   amounts 
to multiplying  the brackets in $\fg(\AA,\OO)$ of two elements of $\fp$
 by $-1$. 
From the description of $\fk$ and $\fp$ we have just given, we see that this twist amounts to 
doing two things. First, twist the Cartan decomposition of $\ft(\OO)$, 
which means that we replace 
$\ft(\OO)=\fso_8$ by its quaternionic form $\fso_{4,4}=\ft(\SO)$. Second, twist the Cayley-Dickson 
process of construction of $\OO$ from $\HH$ by multiplying by $-1$ the product of two elements
in $e\HH$, which amounts to replace $\OO$ by the split Cayley algebra $\SO$. This proves 
the first part of the theorem. The rest of the proof is straightforward.  
\qed 

\subsubsection*{Adjoint varieties and quaternionic symmetric spaces}
Symmetric spaces with quaternionic structures are closely related to adjoint varieties. 
Recall that on a Riemannian manifold $M$, a quaternionic structure is defined as a 
parallel field of quaternion algebras $\HH_x\subset {\rm End}(T_xM)$, $x\in M$, 
such that the unit sphere of $\HH_x$ is contained in the orthogonal group $SO(T_xM)$. 
The dimension of $M$ is then equal to $4m$ for some $m>0$, and its reduced holonomy 
group is contained in $Sp(m)Sp(1)\subset SO(4m)$, where $Sp(m)$ denotes the group
of quaternionic unitary matrices of order $m$ (in particular, $Sp(1)\simeq S^3$ 
is  
the unit sphere in $\HH$). Such matrices act on $\HH^m=\RR^{4m}$ by multiplication 
on the left, which commutes with the scalar multiplication on the right by unitary 
quaternions: the resulting group of orthogonal transformations of $\RR^{4m}$ is 
$Sp(m)Sp(1)$.

This group appears in Berger's classification as one of the 
few possible reduced holonomy groups of nonsymmetric Riemannian manifolds. 
The case of symmetric manifolds was discussed in 
detail by J. Wolf in \cite{wo1}, who proved that there exists exactly one $G$-homogeneous 
symmetric space $M$ with a quaternionic structure for each simple compact Lie group $G$. 
When $G$ is simply connected, we can write $M=G/K.Sp(1)$. Choosing a complex plane in 
$\HH$ amounts to choosing a circle $S^1\subset Sp(1)$, 
and this induces a fibration 
$$X=G/K.S^1\lra M=G/K.Sp(1) $$
whose fibers are two-spheres.
Wolf proved that $X$ is a complex variety, homogeneous under the complexification 
$G^{\CC}$ of the real compact group $G$, and endowed with an invariant contact structure.
Such varieties were previously classified by Boothby, who showed that they are in 
correspondance with simple complex Lie groups   of adjoint type. 
In fact, $X$ is 
  the {\em adjoint variety} $(G^{\CC})^{ad}$, the closed orbit in $\PP\fg^{\CC}$. 

In modern terminology, $X$ is the {\em twistor space} of $M$, and Wolf's work is at the 
origin of this twistor theory, which assigns to each Riemannian manifold $M$ with a
quaternionic structure a complex variety $X$ which is a $S^2$-bundle over $M$, 
defined at each point $x\in M$ by the unit sphere in ${\rm Im} \HH_x$. 
\bigskip

The adjoint varieties are endowed with natural contact 
structures (induced by the Kostant-Kirillov symplectic
structures on the minimal nilpotent coadjoint orbits,  and of which 
the adjoint varieties  are the projectivizations). 
In fact,   any compact Riemannian manifold with exceptional holonomy
$Sp(n)Sp(1)$ and positive scalar curvature   has a twistor
space that is a contact Fano manifold.  
Our understanding of  contact Fano manifolds
has recently improved very much. First of all, there are not
very many:

\begin{theorem} [LeBrun-Salamon \cite{LS}]
Up to biholomorphism, there are only finitely many contact Fano manifolds
of any given dimension.\end{theorem} 

 This lead   to the following conjecture:

\medskip\noindent {\bf Conjecture} (LeBrun-Salamon). {\sl Let $X$ be a smooth complex 
projective contact variety with $b_2(X)=1$. Then $X$ must be the adjoint variety 
of a simple complex Lie group}. 

\medskip
Both the theorem and conjecture can be rephrased purely on the
Riemannian side using the twistor transform.

\medskip
 
In \cite{kebI}, S. Kebekus proved 
that the space 
of minimal rational curves through a fixed 
point of $X$ is a Legendre variety   (except
for the case of $\PP^{2n-1}$), just as is the case for the adjoint varieties. 
 \bigskip

\section{Series of Lie algebras via knot theory and geometry}

\subsection{From knot theory to the universal Lie algebra}
One of the main achievements of the last ten years in topology has been the definition 
of the {\em Kontsevich integral}, which is a universal invariant of finite type for knots. 
The Kontsevich integral associates to each knot a formal series 
in a space of chord diagrams. 
The space is a quotient of the space of
 formal combinations of certain  types of graphs
 by  the AS and IHX relations. From the Kontsevitch integral, 
it is possible to deduce more standard invariants, say numerical or polynomial invariants, 
if one is given a {\em quadratic Lie algebra}, i.e.,
 a finite dimensional Lie algebra endowed 
with an invariant nondegenerate  quadratic form. 
An essential point is that the AS and IHX 
relations can be interpreted respectively as the antisymmetry of the
 Lie bracket  and   the 
Jacobi identity. 

The main reference for the Kontsevich integral is the beautiful article
by Bar-Natan \cite{barnatan}. In it he defines and proves (theorem 4) the mapping
from quadratic Lie algebras equipped with a representation
to a functional on chord diagrams (called a {\it weight system}) 
and makes the conjecture (conjecture 1) that all
  weight systems  come from this mapping. This conjecture
was disproved by P. Vogel in \cite{vog1} which raised the problem
of constructing a category more general than that of quadratic
Lie algebras that would give rise to all weight systems.

In \cite{vog2} Vogel
  defines a candidate for such an object, which he calls the
   {\em universal Lie algebra}.
The universal Lie algebra is a category $\cD$ with the following property: 
for any quadratic Lie (super)algebra $\fg$ over a field $k$, there is a natural functor from 
$\cD$ to the category ${\rm Mod}_k\fg$ of $\fg$-modules. The objects of the category $\cD$ are   finite sets
$[n]=\{1,\ldots ,n\}$, $n\ge 0$. Morphisms in $\cD$ are defined as follows: ${\rm Hom}_{\cD}([p],
[q])$ is the $\ZZ$-module of formal combinations of uni-trivalent
abstract graphs with univalent 
vertices $[p]\cup [q]$, modulo isotopy, and modulo the AS and IHX relations. 
If $\fg$ is a quadratic Lie algebra, a functor from $\cD$ to ${\rm Mod}_k\fg$ is defined at the 
level of objects by sending $[n]$ to $\fg^{\ot n}$. To determine
the 
morphisms, consider some
uni-trivalent graph defining an element of
 ${\rm Hom}_{\cD}([p], [q])$. Possibly after 
isotopy, such a graph can be seen as made of elementary pieces of the form 

$$\begin{array}{ccccc}
\epsfbox{vog.1}\qquad & \epsfbox{vog.2}\qquad & \epsfbox{vog.3}\qquad & 
\epsfbox{vog.4}\quad & \epsfbox{vog.5} \\
bracket\quad  & cobracket\qquad & symmetry\qquad & quadratic\quad  & \quad Casimir \\
 & & & form\qquad & 
\end{array}$$

\noindent
These pieces should be interpreted as indicated, respectively as the Lie bracket 
$\fg^{\ot 2}\ra\fg$, the cobracket  $\fg\ra\fg^{\ot 2}$ 
(the dual to the bracket), 
the symmetry $\fg^{\ot 2}\ra\fg^{\ot 2}$, the quadratic form $\fg^{\ot 2}\ra\fg^{\ot 0}=k$, 
and the Casimir element   $\fg^{\ot 0} \ra \fg^{\ot 2}$ (the dual to the 
quadratic form). Putting together the contributions of 
all its elementary pieces, 
one associates    a morphism from 
$\fg^{\ot p}$ to $\fg^{\ot q}$
to the graph. This morphism   depends only on the isotopy class
of the graph, and   factors through the AS and IHX relations, because of the skew-symmetry 
and Jacobi relation for the Lie bracket in $\fg$. 
\smallskip
Vogel   proves a stronger result   for simple quadratic 
Lie algebras. One can then replace the category $\cD$ by another one $\cD '$, which is 
deduced from $\cD$ by forcing an action of the algebra $\Lambda$ of skew-symmetric connected 
uni-trivalent graphs with exactly three univalent vertices. The functor from $\cD '$
to ${\rm Mod}_k\fg$ is then compatible with the action of 
$\Lambda$ through a character $\chi_\fg :\Lambda \ra k$. 

\smallskip
From the point of view of representation theory, the existence and the properties of the categories 
$\cD$ and $\cD'$ have very interesting consequences. Indeed, each ``decomposition'' of $[p]$ 
in $\cD '$ will imply the existence of a decomposition of $\fg^{\ot p}$ of the same type, {\em for
every simple quadratic Lie algebra}. Here, by a ``decomposition'' of $[p]$, we must understand 
a decomposition of the identity endomorphism of $[p]$ as a sum of idempotents, which will 
correspond to projectors onto submodules of $\fg^{\ot p}$. 

Describing such idempotents is far from easy and   requires quite subtle 
computations with uni-trivalent graphs.  For $p=2$, Vogel obtained the following result.
(Vogel's proof relied on conjecture 3.6 as stated in \cite{vog2}, but  
Vogel himself realized that it is not correct as $\Lambda$, which is conjectured
to be integral, has zero divisors. Nevertheless, a case by case verification  
shows that Vogel's conclusions are indeed correct.)
 
Let $\fg$ be a simple Lie algebra, say over the complex numbers. Then there are decompositions
$$\begin{array}{rcl}
\wedge^2 \fg & = & X_1\op X_2, \\
S^2 \fg & = & X_0\op Y_2\op Y_2'\op Y_2'', 
\end{array}$$
into simple $\fg$-modules (some of which may be zero), with $X_0=\CC$ and $X_1=\fg$. 
Moreover, there exists $\a$, $\b$, $\g$ such that the Casimir operator acts on 
$X_1, X_2, Y_2,  Y_2'$ and $Y_2''$ by multiplication by $t, 2t, 2t-\a, 2t-\b$ and $2t-\g$
respectively, where $t=\a+\b+\g$. 
(In fact the Casimir organizes the representations in
series, e.g., the representations
$X_k$ have Casimir eigenvalue $kt$.)
Finally, the dimensions of these $\fg$-modules are given by 
rational functions in $\a, \b, \g$ (the dimensions of $Y_2'$ and $Y_2''$ are deduced from 
that of $Y_2$ by cyclic permutations of $\a, \b, \g$): 

$${\rm dim} \fg  =  \frac{(\a-2t)(\b-2t)(\g-2t)}{\a\b\g}, $$
$${\rm dim} X_2  =  -\frac{(\a-2t)(\b-2t)(\g-2t)(\a+t)(\b+t)(\g+t)}{\a^2\b^2\g^2}, $$
$${\rm dim} Y_2  =  -\frac{t(\b-2t)(\g-2t)(\b+t)(\g+t)(3\a-2t)}{\a^2\b\g(\a-\b)(\a-\g)}.$$ 

\medskip The scalars $\a, \b, \g$ are readily computed for each simple complex Lie algebra. 
Note that they are defined only up to permutation and multiplication by a same scalar (this
is because the Casimir operator has not been normalized).  We can therefore consider
$(\a,\b,\g)$ as a point in $\PP^2/\cS_3$. 
$$\begin{array}{cccll}
Series & Lie\;algebra & \a & \b & \g \\
 & & & & \\
SL & \fsl_n & -2 & 2 & n\\
OSP & \fso_n, \fsp_{-n} & -2 & 4 & n-4 \\
EXC & \fsl_3 &-2 & 3 & 2\\
&\fg_2 & -3 & 5 & 4\\
&\fso_8 & -2 & 6 & 4\\
 & \ff_4  & -2 & 5 & 6 \\
 & \fe_6  & -2 & 6 & 8 \\
 & \fe_7  & -2 & 8 & 12 \\
 & \fe_8  & -2 & 12  & 20 
\end{array}$$

\smallskip
Note that up in $\PP^2$, the points of the three 
series $SL$, $OSP$ and $EXC$ are located on three
lines, respectively $\a+\b=0$, $2\a+\b=0$, $2\a+2\b-\g=0$. Hence the slogan that {\em there exists
only three complex simple Lie algebras, SL, OSP, and EXC}.  
At the level of representations, it is a classical and well-known fact that the categories of 
modules over the $\fsl_n$, or the $\fso_n$ and $\fsp_{2n}$, have close and precise relationships. 
This is more surprising for the exceptional series.
In \S 4.6 we discuss another collection of Lie algebras that form
a \lq\lq line\rq\rq\ in Vogel's plane.

\begin{remark} In a   spirit similar to Vogel's,
Cvitanovic \cite{cvit}
has proposed a proof of the Killing-Cartan
classification.
 Using  
notation inspired by the Feynman diagrams of quantum field theory,
he represents the invariant tensors of a Lie
algebra by diagrams, and by requring that 
irreducible representations be of integer dimension, he 
deduces severe limits on what simple groups could possibly
exist. His approach dates back to 1976 \cite{C76,C81}. 
The method provides the
correct Killing-Cartan list of all possible simple Lie algebras,
but fails to prove existence.
\end{remark}

\subsection{Vogel's decompositions and Tits correspondences}

As remarked above, Tits correspondences aid one in   decomposing
the tensor powers of nice (e.g., fundamental) representations into
irreducible factors. We illustrate this by showing how to recover
Vogel's decompositions of $S^2\fg$ and $\La^2\fg$. In fact we
can recover Vogel's higher decompositions as well, but we need to
use a slightly more general technique which we call {\it diagram
induction}, see \cite{LMseries}. Using diagram induction, we are also
able to recover the Casimir functions for the representations
occuring in series such as the $X_k$.

If $V$ is a fundamental representation corresponding to a root that
is not short (e.g., the fundamental adjoint representations)
and $X=G/P_{\a}\subset\ppp V$ the closed
orbit, then the Fano variety $\BF_1(X)$ of $\pp 1$'s on $X$ is
$G$-homogeneous by Theorem 2.5, and according to Tits
fibrations,  $\BF_1(X)=G/P_S$
where $S$ is the set of roots joined to $\a$ in $D(\fg )$.
Let $V_2=\langle \BF_1(X)\rangle$ denote
the linear span of the cone over $\BF_1(X)$ in its minimal
homogeneous embedding. Note that
$\BF_1(X)\subset G(2,V)\subset\La^2V$, and we conclude
that $V_2\subset\La^2V$. (To do this, one needs to make sure the
embedding of $\BF_1(X)$ is indeed the minimal one, which is one
reason we need the stronger technique of diagram induction.)

In the case of the adjoint representation, we also know, because
of the cobracket, that $\fg\subset\La^2\fg$. In fact, we
have $\La^2\fg = \fg\op\fg_2$. (In the case of $\fa_n$, $(\fa_n)_2$
is the direct sum of two dual irreducible representations,
 which is also  predicted by Tits fibrations.)
 
 Similiarly, inside $S^2\fg$, we know there is a trivial representation
$X_0$
corresponding to the Killing form, and the Cartan power $\fg\up 2$,
which, breaking Vogel's symmetry, we set $Y_2=\fg\up 2$
so the problem becomes to recover the remaining modules, which we
think of as the \lq\lq primitive\rq\rq\ ones. These are provided in
most cases by Tits fibrations of homogeneous quadrics contained
in $X$. There are up to two such families, we call the larger
one $\fg_Q=Y_2'$ and if there is a second we call it $\fg_{Q'}=Y_2''$,
see \cite{LMseries}. Thus, using Tit's correspondences to find
the primitive factors, one recovers
the Vogel decomposition of $\fg^{\ot 2}$  for
$\fg$ simple (except for $\fc_n$, which
is  tractible using diagram induction).
Even when our methods predict all factors, we would like to emphasize
that they do {\it not} explain why these are the only factors.

  While our observations enable one to recover uniform decompositions
of plethysms in series, the questions of why such series exist and 
how to find them systematically without
looking at the list of simple Lie algebras
 is mysterious using our present techniques.
It would be wonderful and challenging to generalize our perspective sufficiently
to approach Vogel's. One indication that such a generalization
might be possible is that, with the above choices, we obtain
a geometric interpretation of Vogel's parameter $\b$,
namely $\b$ is the dimension of the largest quadric hypersurface
contained in $X_{ad}$ (see \cite{LMseries}). If there is a second
unextendable quadric in $X_{ad}$, then $\gamma$ is its dimension.

\subsection{The exceptional series}
Inspired by Vogel's work,
Deligne \cite{del}, investigated the decompositions
of the tensor powers
of the exceptional Lie algebras into irreducible components, going up
to degree four
(with the help of Cohen, deMan and the computer program LiE, \cite{cdm})
 and giving more explicit
dimension formulas.

\subsubsection{Decomposition formulas}

As with the
Vogel decompositions, these decomposition formulas are all {\em uniform}, 
in the sense, for
example, that the numbers of irreducible components are always the same when the algebra
varies in the exceptional series.  This assertion has to be
 understood with some care: it happens
that some components vanish, or should be taken with a minus sign.
In addition, Deligne and Vogel deal with the algebra twisted
by the symmetry of the marked (for the adjoint representation)
diagram.  
   
   In degree two, one has the decomposition of Vogel, but the
   decomposition of $S^2\fg$ simplifies  because the equation
  $2\a+2\b-\g=0$ implies $Y_2''$ is zero.  In degree three we have
 $$\begin{array}{rcl}
 S^3\fg & = & \fg\op X_2\op A\op Y_3\op Y_3', \\
 \wedge^3\fg & = & \CC\op X_2\op Y_2\op Y_2'\op X_3, \\
 S_{21}\fg&= &2\fg\op X_2\op Y_2 \op Y_2' \op A \op C \op C'
 \end{array}$$
where the (normalized) Casimir operator acts on $X_2, X_3, A, Y_2,
 Y_2', Y_3, Y_3',C,C'$ by 
multiplication by $2, 3, \frac{8}{3}, \frac{2\b+5\g}{3\g},
 \frac{2\a+5\g}{3\g}, 
 2\frac{\b+\g}{\g}, 2\frac{\a+\g}{\g}, \frac{2\b+5\g}{2\g},
 \frac{2\a+5\g}{2\g}$ respectively.
Again, we are able to account for nearly all the factors
and Casimir functions appearing using
diagram induction.  
 
 \smallskip

\subsubsection{Dimension formulas}
  For each exceptional Lie algebra $\fsl_2$, $\fsl_3$, $\fso_8$ 
$\fg_2$, $\ff_4$, $\fe_6$, $\fe_7$, $\fe_8$, let $\l=-3, -2, -1,
-\frac{3}{2}, -\frac{2}{3}, -\frac{1}{2}, -\frac{1}{3},
 -\frac{1}{5}$ respectively. $\lambda$ is a linear function of the length
 of the longest root (with the Casimir normalized to act as the
 identity on $\fg$).
 In
Vogel's parametrization this gives
$(\a,\b,\g)=(\l,1-\l,2)$.

From the decomposition formulas stated above and the knowledge of the Casimir eigenvalue 
on each irreducible component, it is easy (at least with a computer) to calculate 
the dimensions of these components
as  rational functions of $\l$, 
 if one knows the dimension of $\fg$ itself in terms of $\l$. 
   For example,  
$${\rm dim} \fg = -2\frac{(\l+5)(\l-6)}{\l(\l-1)},$$ 
$${\rm dim} X_2 = 5\frac{(\l+3)(\l+5)(\l-4)(\l-6)}{\l^2(\l-1)^2},$$
$${\rm dim} Y_2 = -90\frac{(\l+5)(\l-4)}{\l^2(\l-1)(2\l-1)},$$
$${\rm dim} Y_3 = -10\frac{(\l+5)(5\l-6)(\l-4)(\l-5)(\l-6)}{\l^3(\l-1)^2(2\l-1)(3\l-1)}.$$
A {\em miracle}, Deligne says, is that the
 numerators and denominators of these rational 
functions all factor out into products of linear
 forms with simple integers coefficients! Note that this miracle does
not occur for the Vogel dimension formulas starting with degree three.

In addition,
 there is a striking duality property: the involution $\l\mapsto 1-\l$ takes the dimension 
formula for $Y_2, Y_3$ into those for $Y_2'$, $Y_3'$, while it leaves unchanged the 
dimension formulas for $\fg$, $X_2$, $X_3$, $A$.  

One can attempt to continue the decompositions in degree five,
but there for the first time there appear different representations
with the same Casimir eigenvalue and Deligne's method cannot be
used to compute dimensions. This suggests that if one wants to
continue the formulas, it might be better to restrict to decompositions
into Casimir eigenspaces rather than irreducible modules.

In the last part of \cite{del}, Deligne conjectured the existence of a 
  category,  supposed (just as the universal Lie algebra defined by 
Vogel, but with much stronger properties) to map to the categories of modules over
each exceptional Lie algebra. The properties of this hypothetic category would then
explain at least some of the phenomena discovered by Vogel and Deligne for the 
exceptional series. At the moment, it seems no progress has been made on this 
conjecture. 

\subsubsection{Deligne dimension formulas via triality}
  Deligne's conjecture is elegant while his proof is brute force
computer based. 
In \cite{LMtrial} we   reprove and
extend his formulas using methods whose level of elegance is
somewhere between the conjecture and proof of Deligne. 
The starting
point was our observation  that the parameter $\l$  could be written as 
$\l=-\frac{2}{a+2}$, with $a=-\frac{2}{3}, 1, 2, 4, 8$. (This implies that 
$(\a,\b,\g)=(-2,a+4,2a+4)$.) This indicated
that the relation between exceptional Lie algebras and normed
 division algebras should be 
exploited. The idea was   to use this relation to find a good
description of the exceptional root systems, 
and then apply the Weyl dimension formula. 

The Tits construction is not convenient to describe the exceptional 
root systems and we were led to rediscover  the triality model. 
Recall that for the exceptional series, 
$$\fg(\AA,\OO) = \ft(\AA)\times\ft(\OO)\op (\AA_1\ot\OO_1)\op (\AA_2\ot\OO_2)
\op (\AA_3\ot\OO_3),$$
can be identified with the compact form of $\ff_4, \fe_6, \fe_7, \fe_8$ respectively 
for $\AA = \RR, \CC, \HH, \OO$. In what follows
 we   {\em complexify} this 
construction, without changing  notation. In particular $\ft(\OO)$ is now the 
complex orthogonal Lie algebra $\fso_8(\CC)$. The point is that $\ft(\AA)\times 
\ft(\OO)$ is now a maximal rank reductive subalgebra of $\fg(\AA,\OO)$. In particular,
if we choose Cartan subalgebras $\fh(\AA)\subset\ft(\AA)$ and $\fh(\OO)\subset\ft(\OO)$,
we obtain a Cartan subalgebra $\fh(\AA)\times\fh(\OO)\subset\fg(\AA,\OO)$.
 This provides
a nice description of the root system of $\fg(\AA,\OO)$. There are three kinds of roots:
\begin{itemize} 
\item the {\em heart} of the root system is the set of roots of $\ft(\OO)=\fso_8(\CC)$, 
\item the {\em linear part} is made of the roots of the form $\mu+\nu$, where $\mu$ is a weight 
of $\OO_i$ and $\nu$ a weight of $\AA_i$ for some $i$, 
\item the {\em residue} is the set of roots of $\ft(\AA)$. 
\end{itemize}
When $\AA =\underline 0$,  the root system is of type $D_4$.
When $\AA=\RR$, the root system is of type $F_4$. There is no residue in this case. The 
heart is the set of long roots, and the linear part the set of short roots. Both are root 
systems of type $D_4$. When we change $\AA$ into $\CC, \HH, \OO$, the linear part
{\em expands}, is size being each time multiplied by two, and its roots become long. 
The residue also increases,  but in limited proportions. 

Next  we specify
the positive roots inside the exceptional root systems. They are the roots on which some linear form takes positive values, and we 
choose a linear form that takes very large values on the set $\Delta_+$ of  positive roots 
of $\fso_8(\CC)$. Then the positive roots in the linear part of
 the root system are
the $\mu+\nu$ for $\mu$ belonging to an explicit set $\Sigma$ of weights. This implies that 
the half-sum of the positive roots can be expressed as $\rho = \rho_{\ft(\AA)}+
\rho_{\ft(\OO)}+a\g_{\ft(\OO)}$, 
where $\rho_{\ft(\AA)}\in\fh(\AA)^*$ denotes the half sum of the positive roots in $\ft(\AA)$ and 
$\g_{\ft(\OO)}\in\fh(\OO)^*$ is the sum of the weights belonging to $\Sigma$. 

We are now close to being able to apply the Weyl dimension formula to $\fg(\AA,\OO)$-modules.
What remains to do is to describe the dominant integral weights in $\fh(\AA)^*\times\fh(\OO)^*$. 
At least we can describe the set $C(\OO)\subset\fh(\OO)^*$ of weights that are dominant and 
integral for each $\fg(\AA,\OO)$. Of course,  such weights are be dominant and integral 
in $\fso_8(\CC)$, but this is not   sufficient. 

\begin{prop} The set $C(\OO)\subset\fh(\OO)^*$ is the simplicial cone of nonnegative 
integer linear combinations of the four following weights:
$$
\setlength{\unitlength}{4mm}
\begin{picture}(8.2,2.5)(-5,.3)
\put(-3.7,1){$\o(\fg) =$} 
\put(-.5,1){$1-2$}
\put(1.9,1.5){\line(2,3){.6}}
\put(1.9,1.2){\line(2,-3){.6}}
\put(2.7,2.2){$1$}
\put(2.7,0){$1$}
\end{picture} \quad 
\setlength{\unitlength}{4mm}
\begin{picture}(8.2,2.5)(-5,.3)
\put(-4.5,1){$\o(X_2) =$} 
\put(-.5,1){$2-3$}
\put(1.9,1.5){\line(2,3){.6}}
\put(1.9,1.2){\line(2,-3){.6}}
\put(2.7,2.2){$2$}
\put(2.7,0){$2$}
\end{picture} \quad 
\begin{picture}(8.2,2.5)(-5,.3)
\put(-4.5,1){$\o(X_3) =$} 
\put(-.5,1){$3-4$}
\put(1.9,1.5){\line(2,3){.6}}
\put(1.9,1.2){\line(2,-3){.6}}
\put(2.7,2.2){$2$}
\put(2.7,0){$3$}
\end{picture} \quad 
\begin{picture}(8.2,2.5)(-5,.3)
\put(-4.5,1){$\o(Y_2') =$} 
\put(-.5,1){$2-2$}
\put(1.9,1.5){\line(2,3){.6}}
\put(1.9,1.2){\line(2,-3){.6}}
\put(2.7,2.2){$1$}
\put(2.7,0){$1$}
\end{picture} $$
The lattice in $\ft(\OO)^*$ generated by $C(\OO)$ is not the weight lattice, but the
root lattice of $\fso_8(\CC)$. 
\end{prop}

These weights, that we expressed above in terms of simple roots, 
 occur respectively, $\o(\fg)=\o_2$ as the highest weight of
$\fg$, $\o(X_2)=\o_1+\o_3+\o_4$ as the highest weight of $\La^2\fg$, 
$\o(X_3)=2\o_1+2\o_3$ as the highest weight of $\La^3\fg$, 
and $\o(Y_2')=2\o_1$ as the highest weight of $S^2\fg - \fg\up 2$. 

\smallskip We can now apply the Weyl dimension formula to weights in $C(\OO)$, 
considered for each choice of $\AA$ as integral dominant weights of $\fg(\AA,\OO)$. 
It is essential that the residue of the root system will not contribute. 
We obtain (see \cite{LMtrial} for details):

\begin{theo}\label{super}
 The dimension of the irreducible $\fg(\AA,\OO)$-module
of highest weight $\o\in C(\OO)$ is given by the following formula:
$${\rm dim}\;V_{\o} = \prod_{\a\in\Delta_+\cup\Sigma}
\frac{(a\g_{\ft(\OO)}+\rho_{\ft(\OO)}+\o,\a^{\vee})}
{(a\g_{\ft(\OO)}+\rho_{\ft(\OO)},\a^{\vee})}
\prod_{\beta\in\Sigma}
\frac{\begin{pmatrix}
(a\g_{\ft(\OO)}+\rho_{\ft(\OO)}+\o,\b^{\vee})+\frac{a}{2}-1 \\ (\o,\b^{\vee})
\end{pmatrix}}
{\begin{pmatrix}
(a\g_{\ft(\OO)}+\rho_{\ft(\OO)}+\o,\b^{\vee})-\frac{a}{2} \\ (\o,\b^{\vee})
\end{pmatrix}}.$$
\end{theo}

If $\o= p\o(\fg)+q\o(X_2)+r\o(X_3)+s\o(Y_2^*)$, this formula gives a rational 
function of $a$, whose numerator and denominator are products 
of $6p+12q+16r+10s+24$ linear forms in $a$ with integer coefficients.
Since $\l=-\frac{2}{a+2}$, we obtain formulas of the type of those of Deligne,
and an infinite family of such. 
For example, the $k$-th Cartan power  $\fg(\AA,\OO)^{(k)}$ is $Y_k$ in    
Deligne's notations, and  we get
$$\dim Y_k = \frac{(2k-1)\l-6}{k!\l^k(\l+6)} \prod_{j=1}^k
\frac{((j-1)\l-4)((j-2)\l-5)((j-2)\l-6)}{(j\l-1)((j-1)\l-2)}.$$
Also we can understand why we can expect good dimension formulas for all irreducible 
components of $\fg(\AA,\OO)^{\ot k}$ only for small $k$. This is because when $k$ increases, 
we will unavoidably get components whose highest weight does not come only from $\fso_8(\CC)$, 
but has a contribution from $\fh(\AA)^*$. Then the Weyl dimension formula does not give
a well behaved expression.

While this takes
some of the mystery out
of the dimension formulas,
  the remarkable symmetry $\l\mapsto 1-\l$ noticed
by Deligne, which in our parameter $a$ is $a\mapsto -4\frac{a+3}{a+4}$
remains beyond our understanding.

\subsubsection{Numerology}
In the formula above for $\dim Y_k$, there are very few values of
$\lambda$ for which $\tdim Y_k$ is an integer for all $k$. We are
currently, with B. Westbury \cite{LMW} investigating  
these extra numbers that produce integers and Lie algebras
that go with them. At least in the case $a=6$, one is led to
an algebra of dimension six, the {\it sextonions}, that leads to
a nonreductive row of the magic square between the third and fourth.
Like the {\it odd symplectic groups} of Proctor \cite{proc}
and Gelfand and Zelevinsky \cite{gz}, this series has certain behavior
as if it were a series of reductive Lie algebras.

\subsection{Freudenthal geometries}
The discovery of the Cayley plane $\OO\PP^2$ 
and its automorphism group by Chevalley, Schafer 
and others led to   a period of
 intense activity around the exceptional groups 
and their geometric interpretations in the 1950's 
involving   Freudenthal, Tits and Rozenfeld.

\smallskip
Freudenthal and Rozenfeld   defined explicit varieties 
whose automorphism groups (or groups related to the automorphism
groups) were the exceptional groups.
The starting point is the Jordan 
algebras $\cJ_3(\AA)$.
The projective plane $\AA\PP^2$ is the space
of rank one idempotents of $\cJ_3(\AA)$.
 The subgroup of  $GL(\cJ_3(\AA))$ preserving 
the determinant of order three matrices (which is well defined even 
for $\AA=\OO$) then acts transitively on $\AA\PP^2$. In particular, 
one recovers the transitive action of $E_6$ on 
  $\OO\PP^2$. The subgroup of   $GL(\cJ_3(\AA))$ 
preserving  the determinant and
 the quadratic form $(A,B)={\rm trace} (AB)$ is
  the automorphism group of the Jordan algebra $\cJ_3(\AA)$.
  It  acts irreducibly 
on  $\cJ_3(\AA)_0$, the subspace of trace zero matrices. 
Tits and Freudenthal 
 defined elliptic and projective geometries on $\AA\PP^2$
 with these groups of isometries.  In the elliptic geometry, a {\it point} is the 
same as a {\it line}, and is defined by an element of $\AA\PP^2_0$.   
In the projective geometry,
 a point is defined to be an element of $\AA\PP^2$ and a line
 is determined by
 an element $[\a]$ of the dual projective plane (the rank one idempotents in the
 dual Jordan algebra corresponding to the dual vector space), namely
$\{ X\in \BA\pp 2\mid \a (X)=0\}$.

\smallskip

A synthetic geometry Freudenthal terms {\it symplectic} may be
associated to each of the groups in the third row.
A {\it point} is an element   $[X]\in\PP\fg(\AA,\HH)$,   
  such that ${\rm ad}(X)^2=0$. This condition is
equivalent to requiring   
that $[X]$ is in the adjoint variety $X_{ad}\subset\PP\fg(\AA,\HH)$.
A {\it plane}
   is  is determined by
 an element $[P]\in   \PP\cZ_3(\AA)$, 
  $P\times P=0$. (See \S 4.6 for the definition of the cross product.)
This condition  is equivalent
to saying that $[P]\in G_w(\BA^3,\BA^6)\subset \PP\cZ_3(\AA)$, the
closed orbit. The corresponding plane is  
$\{ [X]\in X_{ad} \mid XP=0\}$.
  A {\it line} is   an intersection of two planes with at least two points.
It is also determined by a point $[\a]\in \BF_1(X)\subset \PP V_2$, where
$V_2$ is the representation defined in \S 4.6  below and $\BF_1(G_w(\BA^3,\BA^6))$ is the closed
orbit (which, as the notation suggests, parametrizes the lines on the
closed orbit $G_w(\BA^3,\BA^6)\subset \PP\cZ_3(\AA)$). 
One can then define incidence rules for these geometric elements which generalize 
the symplectic geometry in $\PP^5$.

\smallskip
For the fourth line of Freudenthal's
magic square Freudenthal defines
 a {\em metasymplectic geometry}. There are now four types of elements, 
points, lines, planes and {\em symplecta}, 
with incidence rules explained
 in detail by Freudenthal. As above, each  geometric element
is a point in the closed 
orbit of the projectivization of a $\fg(\AA,\OO)$-module. 
Here is the table of types of elements:

$$\begin{array}{ccc}
Geometric\;element & Variety & Dimension\\
 & & \\
{\rm Point} & X_{point}\subset\PP \fg_Q & 9a+6\\
{\rm Line} & X_{line}\subset\PP \fg_3 & 11a+9\\
{\rm Plane} & X_{plane}\subset\PP \fg_2 & 9a+11\\
{\rm Symplecta} & X_{symplecta}\subset\PP\fg(\AA,\OO) & 6a+9 
\end{array}$$

\medskip
We see the four generators of the cone $C(\OO)$
in Proposition 4.1 define the four types of 
elements of Freudenthal's metasymplectic geometry! The analogous
results hold for the second and third rows of the magic chart, that
is, Freudenthals elements are the closed orbits in the
projectivizations of the modules whose highest weights
define the fixed cones for our dimension formulas in \cite{LMtrial}.

\subsection{A geometric magic square}
 Consider the following square of complex homogeneous varieties, 
which we investigate in \cite{LMmagic}:

\begin{center}\begin{tabular}{|c|cccc|c} \cline{1-5} 
 & $\RR$ & $\CC$  & $\HH$ & $\OO$ \\ \cline{1-5} 

$\RR$ & $v_2(Q^1)$ & $\PP (T\PP^2)$ & $G_\o(2,6)$ & $\OO\PP^2_{\CC,0}$ &  
section of Severi \\
$\CC$ & $v_2(\PP^2)$ & $\PP^2\times\PP^2$ & $G(2,6)$ & $\OO\PP^2_{\CC}$ & Severi \\ 
$\HH$ & $G_\o(3,6)$ & $G(3,6)$ & $\SS_{12}$ & $E_7/P_7$ & Legendre \\ 
$\OO$ & $F_4^{ad}$ & $E_6^{ad}$ & $E_7^{ad}$ & $E_8^{ad}$ & adjoint \\ \cline{1-5} 
\end{tabular}\end{center}

\medskip
This {\em geometric magic square} is obtained as follows. One begins with the 
adjoint varieties for the exceptional groups: this is the fourth line of the square.
Taking the varieties of lines through a point and applying Theorem 2.5, one obtains 
the third line. The second line is deduced from the third by the same process. 
Then take hyperplane sections to get the first line. We obtain 
projective varieties that are homogeneous under groups  
given by Freudenthal's magic square.

\subsection{The subexceptional  series}

In this final section we consider the series corresponding to
the third row of the extended Freudenthal chart:
$$\fsl_2, \quad \fsl_2\times \fsl_2\times \fsl_2,\quad \fsp_6,  
\quad \fsl_6, \quad \fso_{12}, \quad \fe_7
$$
and view the series from our various perspectives.

\subsubsection{Relation to the universal Lie algebra}
This series also corresponds to a new line in
Vogel's plane, namely $(\a,\b,\g)=(-2,a,a+4)$,
where $a=-\frac 23, 0,1,2,4,8$. Moreover, Vogel's formulas are valid
for the semi-simple Lie algebra  $ \fsl_2\times \fsl_2\times \fsl_2$ despite
the plane being defined only for absolutely simple objects.
Unlike the three lines discovered by Vogel, this line is generic to
order three among actual Lie algebras, in the sense that no
Vogel space is zero except $X_3''$, which is zero for all actual
Lie algebras. We have $\tdim \,\fg (a)=\frac{3(2a+3)(3a+4)}{(a+4)}$.

\subsubsection{Freudenthalia}
For this series there are three preferred representations from
Freudenthal's perspective, corresponding to the ambient spaces
for points, lines, and planes in his incidence geometries.
The points representation is just $\fg=\fg (\BA,\HH)$
and has a nice model thanks to the triality construction.  
For a discussion of the line space, which we will
denote $V_2=V_2(a)$ see \cite{LMmagic}. We have $\tdim V_2(a)= 9(a+1)(2a+3)$.  
The most preferred representation is for the planes, which we
denote $V=V(a)$ and has $\tdim V(a)=6a+8$.  This space  is   the 
complexification of the algebra of {\it Zorn matrices} we encountered in \S 3.8:
$$\cZ_3(\AA) =\Bigg\{\begin{pmatrix} a & X \\ Y & b\end{pmatrix}
,\quad
a,b\in\RR,\; X,Y\in \cJ_3(\AA) \Bigg\}.$$
It can be endowed 
with an algebra structure with multiplication 
$$\begin{pmatrix} a_1 & X_1 \\ Y_1 & b_1\end{pmatrix}
\begin{pmatrix} a_2 & X_2 \\ Y_2 & b_2\end{pmatrix}=
\begin{pmatrix} a_1a_2+\tr(X_1Y_2) & a_1X_2+b_2X_1+Y_1\times Y_2 \\ 
a_2Y_1+b_1Y_2+X_1\times X_2 & b_1b_2+\tr(X_2Y_1)\end{pmatrix}.$$
Here the product $X\times Y$ is defined by the identity ${\rm tr}((X\times Y)Z)=\det (X,Y,Z)$, 
where $\det$ is the polarization of the determinant. 
The algebra is in fact a structurable algebra (see \S 2.5), with
$V(a)=\cZ_3(\AA)=\langle Y_x\rangle \subset T_xX_{adjoint}(\BA,\BO)$
where
$X_{adjoint}(\BA,\BO)$ are the exceptional adjoint varieties.
In fact, Zorn matrices give all simple structurable 
algebras of skew-symmetric rank equal to one, 
and one can show
directly that the derivation algebra of $\cZ_3(\AA)$ is  $\fg(\AA,\HH)$. 

Using this construction one can consider the automorphism group
 of $\cZ_3(\AA)$ in addition to   its derivation 
algebra, and   this can be
done on  any field   (of characteristic 
not equal to $2$ or $3$). 
 Using this, Garibaldi  constructs forms of  $E_7$ over an arbitrary field
\cite{gari}. 
 \bigskip

The closed orbit
$X_{planes}\subset \BP \cZ_3(\AA)$  is
 a natural compactification of the Jordan algebra 
$\cJ_3(\AA)$, and its embedding inside $\PP\cZ_3(\AA)$ 
is given by the translates 
of the determinant on $\cJ_3(\AA)$.
This is the {\em conformal compactification} 
considered by  Faraut and Gindikin in \cite{fg}.
In \cite{LMmagic} we propose a geometric interpretation of
$X_{planes}$ as the Grassmannian of $\BO^3$'s in $\BO^6$
isotropic for a Hermitian symplectic form and use the notation
$G_w(\BO^3,\BO^6)$.

\subsubsection{Decomposition formulas}
The plethysms of $V$ are extraodinarily well behaved in series.
For example, here is a decomposition of $S^kV$ into irreducible components
for each $k$.  
$$\bigoplus_{k\ge 0}t^k S^kV = 
(1-tV)^{-1}(1-t^2\fg)^{-1}(1-t^3V)^{-1}(1-t^4)^{-1}(1-t^4V_2)^{-1}.$$
The right hand side in the formula is to be
expanded out in geometric
series and
multiplication of
representations is taken
in terms of Cartan products $V_{\lambda}V_{\mu}=V_{\lambda +\mu}$.
The case of $\fsl_2\times\fsl_2\times\fsl_2$ is somewhat special, since the formula 
above holds true only if the natural action of $\cS_3$ is taken into account.  
See \cite{LMseries}.

Each of the factors in the rational function above can be accounted for in terms
of diagram induction. In particular, $\fg=V_Q$, the ambient space
for the variety parametrizing the $G$-homogeneous quadrics on $X_{planes}$
and $V_2=\langle \BF_1(X_{planes})\rangle$ is the ambient space for
the variety parametrizing the lines on $X_{planes}$.

\subsubsection{Dimension formulas}
Our method for applying the Weyl dimension formula in series
also works for the three distinguished representations.
For example, we have
$$
\dim V^{(k)}   =   \frac{2a+2k+2}{a+1}\frac{\binom{k+2a+1}{2a+1}
\binom{k+\frac{3a}{2}+1}{\frac{3a}{2}+1}}{\binom{k+
\frac{a}{2}+1}{\frac{a}{2}+1}},
$$
where the binomial coefficients
are defined by
$\binom{k+x}{k} = (1+x)\cdots (k+x)/k!$  and thus are rational polynomials 
of degre $k$ in $x$.  

Note that since $V\up k$ is the complement of $I_k(X_{planes})$ in
$S^kV^*$, the above formula also gives the Hilbert functions of
the varieties $X_{planes}$ in a uniform manner.  

\subsubsection{Zakology}
The varieties $X_{planes}\subset \BP V$ are important
for the following geometric
classification problem: Given a smooth variety
$X\subset \BP V$, one defines its dual variety $X^*\subset\BP V^*$
to be the set of hyperplanes tangent to $X$. Usually the degree
of $X^*$ is quite large with respect to $X$
and Zak has proposed the problem of classifying smooth varieties
whose duals are of low degree \cite{zakdual}.

Here the duals of $X_{planes}\subset \BP V$ are of degree four.
Moreover, they are the tangential varieties of the closed orbits
in the dual projective space, i.e., we have $X^*_{planes}\simeq \tau (X_{planes})$.

The first variety in the series is $X_{planes}(-2/3)=v_3(\pp 1)$ and the equation
of its dual is the classical discriminant of a cubic polynomial.
Using $\BA$-valued variables, we   write the equations for the 
duals in a uniform fashion, see \cite{LMmagic}. In particular, we
characterize the other quartics by their restriction to the
preferred subspace. The second variety in the series is
$X=Seg(\pp 1\times \pp 1\times \pp 1)$ and its dual is Cayley's
{\it hyperdeterminant},
 see \cite{gkz}. Our restriction result gives a new characterization
of the hyperdeterminant.

\subsubsection{Orbits}
The orbit structure is also uniform. 
\begin{prop} \cite{LMmagic} For each of the 
varieties $X_{planes}=G_w(\BA^3, \BA^6) \subset \PP V $ there are
exactly four
orbits, the closures of which are ordered by inclusion:
$$
G_w(\BA^3, \BA^6)\subset\s_+(G_w(\BA^3, \BA^6))\subset\tau  (G_w(\BA^3,
\BA^6))\subset \PP V.$$ 
The dimensions are respectively
$3a+3, 5a +3$ and $6a+6$.
\end{prop}
Needless to say,
 the singular orbit closures have uniform desingularizations
by Kempf's method.
  
\medskip
Once again, the triality model plays a   unifying role for geometric 
and representation geometric phenomena which at first glance  seem   
sporadic. Another striking example of this role, which is currently under 
investigation, concerns nilpotent orbits in exceptional Lie algebras and the 
associated unipotent characters of exceptional Chevalley groups \cite{LMW2}.

\vspace{2cm}

\begin{tabular}{lll}
Joseph M. Landsberg & \hspace*{1cm} & Laurent Manivel \\
School of Mathematics, & & Institut Fourier, UMR 5582 du CNRS \\
Georgia Institute of Technology, & & Universit\'e Grenoble I, BP 74 \\
Atlanta, GA 30332-0160 & & 38402 Saint Martin d'H\`eres cedex \\
USA & & FRANCE \\
{\rm E-mail}: jml@math.gatech.edu & &
{\rm E-mail}: Laurent.Manivel@ujf-grenoble.fr 
\end{tabular}

\end{document}